\documentclass[12pt]{article}
\usepackage{amsmath}
\usepackage{latexsym}
\usepackage{amssymb}
\newtheorem{thm}{Theorem}[section]
\newtheorem{la}[thm]{Lemma}
\newtheorem{Defn}[thm]{Definition}
\newtheorem{Remark}[thm]{Remark}
\newtheorem{prop}[thm]{Proposition}
\newtheorem{cor}[thm]{Corollary}
\newtheorem{Example}[thm]{Example}
\newtheorem{Number}[thm]{\!\!}
\newenvironment{defn}{\begin{Defn}\rm}{\end{Defn}}

\newenvironment{rem}{\begin{Remark}\rm}{\end{Remark}}
\newenvironment{numba}{\begin{Number}\rm}{\end{Number}}
\newenvironment{proof}{{\noindent\bf Proof.}}%
                  {\nopagebreak\hspace*{\fill}$\Box$\medskip\medskip\par}   
\newcommand{\Punkt}{\nopagebreak\hspace*{\fill}$\Box$}
\newcommand{\wb}{\overline}

\newcommand{\at}{\symbol{'100}}

\newcommand{\tensor}{\otimes}

\newcommand{\mto}{\mapsto}
\newcommand{\isom}{\cong}
\newcommand{\N}{{\mathbb N}}

\newcommand{\K}{{\mathbb K}}
\newcommand{\Q}{{\mathbb Q}}
\newcommand{\Z}{{\mathbb Z}}

\newcommand{\cp}{{\mathfrak p}}
\DeclareMathOperator{\Aut}{Aut}

\newcommand{\one}{{\bf 1}}
\newcommand{\sub}{\subseteq}
\DeclareMathOperator{\GL}{GL}
\DeclareMathOperator{\im}{im}

\DeclareMathOperator{\Spann}{span}

\newcommand{\sbull}{{\scriptscriptstyle \bullet}}
\DeclareMathOperator{\Int}{Int}
\DeclareMathOperator{\tor}{tor}
\DeclareMathOperator{\dv}{div}

\newcommand{\semid}{\mbox{$\times\!$\rule{.15 mm}{2.07 mm}}}
\begin{document}
\begin{center}
{\Large\bf Classification of the simple factors appearing
in composition series of totally
disconnected\\[2.1mm]
contraction groups}\\[6mm]
\renewcommand{\thefootnote}{\fnsymbol{footnote}}
{\bf Helge Gl\"{o}ckner and
George A. Willis%
\footnote{Research
supported by DFG grant 447 AUS-113/22/0-1
and ARC grants LX 0349209 and DP0208137}}\vspace{4mm}
\end{center}
\renewcommand{\thefootnote}{\arabic{footnote}}
\setcounter{footnote}{0}
\begin{abstract}\vspace{1mm}
\hspace*{-7.2 mm}
Let $G$ be a totally disconnected, locally compact
group admitting a contractive automorphism~$\alpha$.
We prove a Jordan-H\"{o}lder
theorem for
series of $\alpha$-stable
closed subgroups of~$G$,
classify all possible
compo\-sition factors and deduce consequences for the
structure of~$G$.\vspace{3mm}
\end{abstract}
{\footnotesize {\em Classification}:
22D05 (primary); %gen props and structure of lcp gps
20E15,           %chains and lattices of subgroups, subnormal subgroups
20E36\\[2.5mm]   %gen thms concerning automorphisms
{\em Key words}:
totally disconnected group,
automorphism, contraction group,
normal series, Jordan-H\"{o}lder theorem,
composition factor, simple group, classification,
$p$-adic Lie group,
shift, restricted product, pro-discrete group,
torsion group, divisibility}\vspace{5mm}
\begin{center}
{\bf\Large Introduction}
\end{center}
A \emph{contraction group} is a pair $(G,\alpha)$,
where $G$ is a topological group
and $\alpha\colon G\to G$ a contractive
automorphism, meaning that
$\alpha^n(x)\to 1$ as $n\to \infty$,
for each $x\in G$. Contraction groups
arise in probability theory on
locally compact groups (see, e.g., \cite{HaS}),
representation theory (\cite{Mar}, \cite{Mau},
\cite{Moo}, \cite{Wan}),
and the structure theory of locally
compact groups initiated in \cite{Wi1} (see \cite{BaW}
and \cite{GLA}).
It is known from the work of Siebert
that every locally compact contraction group
is a direct product $G=G_e\times D$
of a connected group~$G_e$
and an $\alpha$-stable totally disconnected group~$D$,
whence the study of
locally compact contraction groups
splits into the two extreme
cases of connected groups
and totally disconnected groups
(see \cite[Proposition~4.2]{Sie}).
Siebert characterized the connected
locally compact
contraction groups;
they are, in particular, simply connected,
nilpotent real Lie groups.
He also provided some basic
information concerning the totally disconnected case.
In the present article, we complete
the picture by discussing the fine structure
of a totally disconnected, locally compact
contraction group~$G$.
We show the existence of a composition
series
\[
\one=G_0 \,\triangleleft\; \cdots \;  \triangleleft \, G_n=G
\]
of $\alpha$-stable closed subgroups
of~$G$,
prove a Jordan-H\"{o}lder Theorem
for the topological factor groups
and find the
possible
composition factors $G_j/G_{j-1}$.
Any such is a \emph{simple}
contraction group in the sense
that it does not have a non-trivial,
proper, $\alpha$-stable
closed normal subgroup.
Our first main result is a
classification
of the simple, totally disconnected
contraction groups.\\[3mm]
{\bf Theorem A.}
\emph{Let $(G,\alpha)$ be a simple,
totally disconnected, locally compact
contraction group.
Then $G$ is a torsion group or
torsion-free, and $(G,\alpha)$ is of the following form:}
\begin{itemize}
\item[(a)]
\emph{If $\,G$ is a torsion group,
then $(G,\alpha)$ is isomorphic
to a restricted product $F^{({-\N})}\times F^{\N_0}$
with the right shift,
for a finite simple group~$F$.}
\item[(b)]
\emph{If $\,G$ is torsion-free,
then $(G,\alpha)$ is isomorphic
to a finite-dimensional
$p$-adic vector space,
together with a contractive linear
automorphism
which does not leave any non-trivial, proper
vector subspace invariant.}
\end{itemize}
\emph{Conversely, all of the contraction groups
described in {\rm (a)} and {\rm (b)} are simple.}\\[2.5mm] 
We remark that
the contractive linear automorphisms occurring in (b)
can be characterized
in terms of their rational normal form
(cf.\ Proposition~\ref{ratform}).
The classification has important consequences for
general contraction groups.
Our second main result is the following structure~theorem.\\[3mm]
{\bf Theorem B.}
\emph{Let $(G,\alpha)$ be a totally disconnected,
locally compact contraction group.
Then the set $\tor(G)$ of torsion elements
and the set $\dv(G)$ of infinitely divisible
elements are
$\alpha$-stable
closed subgroups of~$G$, and}
\[
G\;=\; \tor(G)\times\dv(G)
\]
\emph{$($internally$)$ as a topological group.
Furthermore, $\dv(G)$ is a direct
product of $\alpha$-stable, nilpotent
$p$-adic Lie groups $G_p$ for certain primes~$p$},
\[
\dv(G)\;=\; G_{p_1}\times\, \cdots\, \times G_{p_n}\,.
\]
Each $G_p$ actually is the group
of $\Q_p$-rational points of a unipotent
linear\linebreak
algebraic group defined over $\Q_p$,
by \cite[Theorem~3.5\,(ii)]{Wan}.\\[2.5mm]
{\bf Organization of the article.}
Sections~\ref{secprel}
and~\ref{sectopop} are of a preparatory nature.
In Section~\ref{secprel},
we compile
several basic facts concerning contraction
groups.\linebreak
In Section~\ref{sectopop},
we fix our terminology concerning
topological groups with operators
and formulate a criterion
for the validity of a Jordan-H\"{o}lder Theorem,
which can be verified
for the cases of relevance (in Section~\ref{secjhold}).
This is quite remarkable,
because composition series
can rarely be used with profit
in the theory of topological groups
(typically they need not exist,
and if they do, then uniqueness
of the composition factors
cannot be insured).
Sections~\ref{secprodis} and~\ref{sectools}
prepare the proof
of the classification
(given in Section~\ref{seccompfac}).
Notably, we show
there that every simple totally
disconnected contraction group
is \emph{pro-discrete},
i.e.,
every identity neighbourhood
contains an open normal subgroup.
As a tool for the proof
of Theorem~B, we explain in Section~\ref{sec:series}
how a canonical series
of $\alpha$-stable normal subgroups
can be associated with a contraction group.
The proof of the Structure Theorem
is outlined in Section~\ref{sec:outline}
and details are provided
in Sections~\ref{secextprod}--\ref{sec:divisible}.\\[4mm]
\emph{Acknowledgement.} We are grateful to Udo Baumgartner
for discussions in the early stages of this work.

\section{Preliminaries}\label{secprel}
Let us agree on the following conventions
concerning subgroups and automorphisms of a topological
group~$G$.
All automorphisms $\alpha$ of~$G$
are assumed bicontinuous.
A subgroup $H\leq G$
is called \emph{$\alpha$-stable} if $\alpha(H)=H$
while it is \emph{$\alpha$-invariant}
if $\alpha(H)\sub H$.
If $H$ is $\alpha$-stable
for each automorphism~$\alpha$
of~$G$, then $H$ is called
\emph{topologically characteristic}.
It is \emph{topologically fully invariant}
if it is invariant under each
endomorphism of the topological group~$G$.
An automorphism $\alpha$ of~$G$
is called \emph{contractive}
if $\alpha^n(x)\to 1$ for all $x\in G$.
The module of an automorphism~$\alpha$ of a locally compact
group~$G$ is defined as $\Delta_G(\alpha):=
\frac{\lambda(\alpha(U))}{\lambda(U)}$,
where $U\sub G$ is any non-empty, relatively compact, open
subset and $\lambda$ a left invariant Haar measure
on~$G$.\\[2.5mm]
We recall various important facts concerning
contractive automorphisms.
%
%\ma{basefacts}
\begin{prop}\label{basefacts}
For each totally disconnected, locally compact
group~$G$ and contractive automorphism
$\alpha\colon G\to G$,
the following holds:
\begin{itemize}
\item[\rm (a)]
$\alpha$ is compactly contractive, i.e.,
for each compact subset $K\sub G$ and
identity neighbourhood $W\sub G$,
there is $N\in\N$ such that
$\alpha^n(K)\sub W$ for all $n\geq N$
$($and hence also $K\sub \alpha^{-n}(W))$.
\item[\rm (b)]
$G$ has a compact, open subgroup $W$
such that $\alpha(W)\sub W$.
If $G$ is pro-discrete,
then $W$ can be chosen as a normal subgroup of~$G$.
\item[\rm (c)]
If $G\not=\one$, then $G$ is neither discrete
nor compact.
\item[\rm (d)]
If $W\sub G$ is a relatively compact, open
identity neighbourhood, then
$\{\alpha^n(W)\colon n\in \N_0\}$
is a basis of identity neighbourhoods,
and $G=\bigcup_{n\in \N_0}\alpha^{-n}(W)$.
In particular, $G$ is
first countable and $\sigma$-compact.
\item[\rm (e)]
The module $\Delta_G(\alpha^{-1})$
is an integer $\geq 2$.
\item[\rm (f)]
If $G$ has a compact, open, normal subgroup,
then $G$ is pro-discrete.
\end{itemize}
\end{prop}
\begin{proof}
(a) See \cite[Proposition~2.1]{Wan} or \cite[Lemma~1.4\,(iv)]{Sie}.

(b) See \cite[Lemma~3.2\,(i) and Remark~3.4\,(2)]{Sie}.

(c) See \cite[3.1]{Sie}.

(d) follows directly from (a); see also
\cite[1.8\,(a)]{Sie}.

(e) For $W$ as in (b), $W$ is an open subgroup
of the compact group $\alpha^{-1}(W)$ and thus
\[
\Delta_G(\alpha^{-1})\; =\;
\frac{\lambda(\alpha^{-1}(W))}{\lambda(W)}
\; =\; [\alpha^{-1}(W):W]
\; \in \; \N \, .
\]
If $[\alpha^{-1}(W):W]$ was~$1$, then we would have
$W=\alpha^{-1}(W)$ and thus
$G=\bigcup_{n\in \N_0}\alpha^{-n}(W)=W$,
whence $G$ would be compact,
contradicting~(c).

(f) See \cite[Remark~3.4\,(2)]{Sie}.
\end{proof}
If $G$ is a topological group
and $\alpha\colon G\to G$ a contractive automorphism,
it is convenient to call
$(G,\alpha)$ a \emph{contraction group}.\\[2.5mm]
{\footnotesize\bf In the following, all contraction groups
are assumed locally compact
and totally disconnected,
unless the contrary is stated.}\\[2.5mm]
The proof of the classification
hinges on the theory of analytic pro-$p$-groups.
We refer to~\cite{DMS} for
background information.
Generalities concerning $p$-adic Lie groups
can also be found in \cite{Bo2} and~\cite{Ser}.
Standard facts from the theory of
pro-finite groups and their Sylow subgroups
(as provided in \cite{RaZ} or \cite{Wsn})
will be used freely.
We shall say that a topological group~$G$
is \emph{locally pro-$p$}
if it has an open subgroup
which is a pro-$p$-group.
\section{Topological groups with operators}\label{sectopop}
We are interested
in series of $\alpha$-stable
closed subgroups of contraction groups, but also
in series of $\alpha$-stable, closed,
normal subgroups.
Topological groups
with operators provide
the appropriate
language to deal with both cases simultaneously.
They also enable us to formulate
sufficient conditions for the validity of
a Jordan-H\"{o}lder Theorem.
\begin{defn}
Let $\Omega$ be a set.
A \emph{topological $\Omega$-group}
is a Hausdorff topological group~$G$,
together with a map $\kappa\colon
\Omega\times G\to G$, $(\omega,g)\mto \kappa(\omega,g)=:\omega.g$
such that $\kappa(\omega,\sbull)\colon G\to G$
is a continuous endomorphism of~$G$, for each
$\omega\in \Omega$.
A subgroup $H$ of a topological $\Omega$-group~$G$
is called an \emph{$\Omega$-subgroup}
if it is closed and $\Omega.H\sub H$.
A continuous homomorphism
$\phi\colon G\to H$ between topological
$\Omega$-groups
is called an \emph{$\Omega$-homomorphism}
if $\phi(\omega.g)=\omega.\phi(g)$ for all
$\omega\in \Omega$ and $g\in G$.
\end{defn}
Each $\Omega$-subgroup $H\leq G$ of a topological
$\Omega$-group~$G$
and each quotient $G/N$ by a normal
$\Omega$-subgroup is a topological
$\Omega$-group in a natural way.
\begin{rem}
If $(G,\alpha)$ is a totally disconnected
contraction group, we shall turn $G$
into a topological $\Omega$-group
in two ways:
\begin{itemize}
\item[(a)]
$\Omega=\langle \alpha\rangle\leq \Aut(G)$.
Then $\langle\alpha\rangle$-subgroups
are closed $\alpha$-stable subgroups.
\item[(b)]
$\Omega=\langle \Int(G)\cup\{\alpha\}\rangle\leq \Aut(G)$,
where $\Int(G)$ is the group
of inner automorphisms of~$G$.
In this case, the $\Omega$-subgroups of~$G$
are the closed $\alpha$-stable normal subgroups of~$G$.
\end{itemize}
\end{rem}
Let $G$ be a topological $\Omega$-group. As expected,
a series
%
%\ma{formseries}
\begin{equation}\label{formseries}
\one\;=\; G_0\; \triangleleft \; G_1\;\triangleleft \; \cdots\;\triangleleft  \;G_n=G
\end{equation}
is called an \emph{$\Omega$-series}
if each $G_j$ is an $\Omega$-subgroup
of~$G$ (and hence closed).
An $\Omega$-composition
series is an $\Omega$-series (\ref{formseries})
without repetitions
which does not admit a proper refinement.
Two $\Omega$-series ${\bf S}$ and ${\bf T}$
are $\Omega$-isomorphic
if there is a bijection from the set of
factors of ${\bf S}$ onto the set of factors
of~${\bf T}$ such that corresponding factors
are $\Omega$-isomorphic as topological
$\Omega$-groups.
\begin{defn}
A totally disconnected contraction group
$(G,\alpha)$ is called \emph{simple}
if it is a simple topological $\langle\alpha\rangle$-group,
that is,
$G\not=\one$ and $G$ has no
$\alpha$-stable closed normal subgroups except for $\one$ and $G$.
\end{defn}
Evidently,
an $\langle\alpha\rangle$-series
$\one=G_0\,\triangleleft\,\cdots\, \triangleleft \, G_n=G$ of a contraction group
$(G,\alpha)$ is an
$\langle\alpha\rangle$-composition series
if and only if all factors $G_j/G_{j-1}$
are simple contraction groups.\\[2.5mm]
We now formulate a criterion for the validity
of the Jordan-H\"{o}lder Theorem.
%
%
%
%\ma{preJordH}
\begin{prop}\label{preJordH}
Let $G$ be a $\sigma$-compact,
locally compact group with set of operators~$\Omega$.
Assume that $G$ satisfies the following
``closed product property'':
For all $\Omega$-subgroups $H_1, H_2\leq G$
such that $H_2$ normalizes~$H_1$,
the product $H_1H_2$ is closed in~$G$
$($and hence an $\Omega$-subgroup$)$.
Then the following holds:
\begin{itemize}
\item[\rm (a)]
$($Schreier Refinement Theorem$)$
Any two $\Omega$-series
of $G$ have $\Omega$-isomor\-phic refinements.
\item[\rm (b)]
$($Jordan-H\"{o}lder Theorem$)$
If ${\bf S}$ is an $\Omega$-composition series
and ${\bf T}$ is any $\Omega$-series
of~$G$, then ${\bf T}$ has a refinement which
is an $\Omega$-composition series
and is $\Omega$-isomorphic
with ${\bf S}$.
In particular, any two $\Omega$-composition series
of~$G$ are $\Omega$-isomorphic.
\end{itemize}
\end{prop}
\begin{proof}
If
$S\leq G$,
$H \leq S$ and $N\triangleleft S$
are closed subgroups
such that $HN$ is closed,
then $H$ and $HN/N$
are $\sigma$-compact,
locally compact groups.
Therefore $H/(H\cap N)\to HN/N$,
$x(H\cap N)\mto xN$
is a topological isomorphism,
by \cite[(5.33)]{HaR}.
Using this information and
the closed product
property, we find that all
subgroups occurring in the standard
proofs of the Zassenhaus Lemma and the Schreier
Refinement Theorem
(as in \cite[3.1.1--3.1.2]{Rob})
are closed and that all relevant abstract isomorphisms
are isomorphisms of topological groups.
Thus (a) holds. Part\,(b) is a direct consequence.
\end{proof}
\section{\hspace*{-2.6mm}The\! Jordan-H\"{o}lder\! Theorem\! for\!
contraction\protect\linebreak \hspace*{-2.4mm}groups}\label{secjhold}
We now verify the closed product
property (from Proposition~\ref{preJordH})
for $\alpha$-stable closed subgroups
of a totally disconnected contraction group $(G,\alpha)$.
As a consequence, a Jordan-H\"{o}lder Theorem
holds for contraction groups.
We shall also see that $\langle\alpha\rangle$-composition series
always exist.\\[2.5mm]
The following proposition is
one of the main technical tools of this article.
It ensures that an $\alpha$-invariant closed
subgroup of a contraction group $(G,\alpha)$
always is an open subgroup of a suitable
$\alpha$-stable closed subgroup of~$G$.
%
%
%\ma{inflate}
\begin{prop}\label{inflate}
Let $(G,\alpha)$ be a totally disconnected, locally compact
contraction group
and $H\leq G$ be a closed
subgroup such that $\alpha(H)\sub H$.
Then
\begin{itemize}
\item[\rm (a)]
$S := \bigcup_{n\in \N_0}\alpha^{-n}(H)$
is a closed $\alpha$-stable
subgroup of~$G$,
and $H$ is open in~$S$.
Furthermore:
\item[\rm (b)]
There exists a compact, open subgroup
$K\leq H$ such that $\alpha(K)\sub K$.
\item[\rm (c)]
If $H$ is compact, then $[H:\alpha(H)]<\infty$.
\item[\rm (d)]
If $H$ is normal in $G$, then also $S$ is normal
in~$G$.
\end{itemize}
\end{prop}
\begin{proof}
(a) and (d):
Since $\alpha(H)\sub H$, we have
$\alpha^{-j}(H)\sub \alpha^{-j-1}(H)$
for all $j\in \N_0$,
entailing that $S$
is a subgroup of~$G$
(which is normal if so is~$H$).
By Proposition~\ref{basefacts}\,(b),
there exists a compact, open subgroup
$W\leq G$ such that $\alpha(W)\sub W$.
We claim that
%
%\ma{essnz}
\begin{equation}\label{essnz}
S\cap \alpha^{n_0}(W)\;=\;H \cap\alpha^{n_0}(W)
\end{equation}
for some $n_0\in\N_0$.
If this claim is true, then
$S\cap \alpha^{n_0}(W)$ is a compact
identity neighbourhood in~$S$,
and hence $S$ is locally compact
and thus closed in~$G$.
Furthermore, $H$ is open in $S$,
as it contains the open set
$S\cap \alpha^{n_0}(W)$.
To prove (a), it therefore
only remains to establish~(\ref{essnz}).
We proceed in steps.\vspace{.5mm}

{\bf Step~1.}
We first note that the indices
$\ell_n:=[H\cap\alpha^n(W):H\cap \alpha^{n+1}(W)]
=[(H\cap\alpha^n(W))\alpha^{n+1}(W):\alpha^{n+1}(W)]
\leq [\alpha^n(W):\alpha^{n+1}(W)]=[W:\alpha(W)]$
are bounded. Furthermore,
the sequence $(\ell_n)_{n\in \N_0}$ is monotonically
increasing, as
\begin{eqnarray*}
[(\!H\cap \alpha^n(W))\alpha^{n+1}(W):\alpha^{n+1}(W)]
&\!\!\!\!=\!\!\!\!&
[(\alpha(H)\cap \alpha^{n+1}(W))\alpha^{n+2}(W)\!:\!\alpha^{n+2}(W)]\\
&\!\!\!\!\leq\!\!\!\!&
[(H \cap \alpha^{n+1}(W))\alpha^{n+2}(W):\alpha^{n+2}(W)].
\end{eqnarray*}
Here, we applied the automorphism $\alpha$ to all subgroups,
and then used that $\alpha(H)\sub H$.
As a consequence, $(\ell_n)_{n\in \N}$
becomes stationary; there are
$n_0\in \N_0$ and $\ell\in \N$ such that
$[H\cap\alpha^n(W):H\cap \alpha^{n+1}(W)]=\ell$
for all $n\geq n_0$.\vspace{.5mm}

{\bf Step~2.}
\emph{$\alpha^m(H\cap\alpha^n(W))=H\cap \alpha^{n+m}(W)$,
for all $n\geq n_0$ and $m\in \N_0$.}
To see this, note first
that $\alpha^m(H\cap \alpha^n(W))
=\alpha^m(H)\cap \alpha^{n+m}(W)
\sub H\cap \alpha^{n+m}(W)$.
Since $[(H\cap \alpha^n(W))\alpha^{n+1}(W):\alpha^{n+1}(W)]=\ell$,
we find
$x_1,\ldots,x_\ell\in H\cap\alpha^n(W)$
such that
$(H\cap \alpha^n(W))\alpha^{n+1}(W)/\alpha^{n+1}(W)
=\{x_1\alpha^{n+1}(W),\ldots, x_\ell\alpha^{n+1}(W)\}$
and thus
$x_i^{-1}x_j\not\in \alpha^{n+1}(W)$
whenever $i\not= j$.
For each $k\in \N_0$, we then have
\[
\alpha^k(x_1),\ldots,\alpha^k(x_\ell)\,\in\,
\alpha^k(H)\cap\alpha^{n+k}(H)\,\sub\,
H\cap\alpha^{n+k}(W)
\]
and $\alpha^k(x_i)^{-1}\alpha^k(x_j)\not\in\alpha^{n+k+1}(W)$
if $i\not= j$,
entailing that
\[
(H\cap \alpha^{n+k}(W))\alpha^{n+k+1}(W)/\alpha^{n+k+1}(W)
=
\{\alpha^k(x_i)\alpha^{n+k+1}(W):
i=1,\ldots, \ell\}
\]
for all $k\in \N_0$.
Thus, given $y\in H\cap \alpha^{n+m}(W)$,
we find $i_0\in \{1,\ldots,\ell\}$
such that $\alpha^m(x_{i_0})^{-1}y\in H\cap \alpha^{n+m+1}(W)$.
Next, we find $i_1\in \{1,\ldots,\ell\}$
such that $\alpha^{m+1}(x_{i_1})^{-1}\alpha^m(x_{i_0})^{-1}y
\in H\cap \alpha^{n+m+2}(W)$.
Proceeding in this way, we obtain
a sequence $(i_k)_{k\in \N_0}$ in
$\{1,\ldots,\ell\}$
such that
\[
\alpha^{m+k}(x_{i_k})^{-1}\cdots
\alpha^m(x_{i_0})^{-1}y\in H\cap \alpha^{n+m+k+1}(W)
\]
for all $k\in \N_0$ and thus
$\alpha^m(x_{i_0}\alpha(x_{i_1})\cdots \alpha^k(x_{i_k}))
\in y \alpha^{n+m+k+1}(W)$,
whence
\[
y\,=\,\lim_{k\to\infty}\alpha^m(x_{i_0}\alpha(x_{i_1})\cdots \alpha^k(x_{i_k}))
\,=\,\alpha^m(x)
\]
with $x:=\lim_{k\to\infty}
x_{i_0}\alpha(x_{i_1})\cdots \alpha^k(x_{i_k})\in H\cap\alpha^n(W)$
(as $H$ is closed).\vspace{.5mm}

{\bf Step~3.}
By Step~2, we have $\alpha^{-m}(H)\cap \alpha^{n_0}(W)
=H\cap \alpha^{n_0}(W)$ for all $m\in \N_0$,
whence
$S\cap \alpha^{n_0}(W)=H\cap\alpha^{n_0}(W)$.
Thus (\ref{essnz})
(and thus (a))~hold.

(b) Let $V\sub S$ be a compact, open subgroup
such that $\alpha(V)\sub V$.
Since $H$ is open in~$V$,
we have $K:=\alpha^n(V)\sub H$
for some $n\in \N$, and this is a subgroup
with the desired properties.

(c) Since $\alpha|_S$ is an automorphism
of~$S$ and $H$ is open in~$S$,
the image $\alpha(H)$ is open in~$H$
and therefore has finite index if~$H$ is compact.
\end{proof}
%
%
%\ma{havecpp}
\begin{cor}\label{havecpp}
Let $(G,\alpha)$ be a totally disconnected
contraction group. Then the following holds:
\begin{itemize}
\item[\rm (a)]
$H_1H_2$ is closed in~$G$,
for any closed
subgroups $H_1, H_2\leq G$ such that
$\alpha(H_1)\sub H_1$, $\alpha(H_2)\sub H_2$,
and $H_2$ normalizes~$H_1$.
\item[\rm(b)]
For both $\Omega=\langle\alpha \rangle\leq \Aut(G)$
and $\Omega=\langle \Int(G)\cup\{\alpha\}\rangle$,
the topological $\Omega$-group~$G$
has the closed product property.
\end{itemize}
\end{cor}
\begin{proof}
(a) By Proposition~\ref{inflate}\,(b),
there exists a compact, open subgroup
$W\sub H_2$
such that $\alpha(W)\sub W$.
Then $W$ normalizes~$H_1$, and thus
$H:=H_1W$ is a subgroup of~$G$.
We have $\alpha(H_1W)=\alpha(H_1)\alpha(W)\sub H_1W$,
and furthermore $H_1W$
is closed in~$G$ because $H_1$ is closed
and $W$ is compact~\cite[(4.4)]{HaR}.
Hence $S:=\bigcup_{n\in \N_0}\alpha^{-n}(H)$
is closed in~$G$ and~$H$ is
open in~$S$, by
Proposition~\ref{inflate}\,(a).
Since $H_1\sub \alpha^{-1}(H_1)$
and $H_2\sub \bigcup_{n\in \N_0}\alpha^{-n}(W)$,
we have $H_1H_2\sub S$.
As $H\sub H_1H_2$,
the product $H_1H_2$
is an open subgroup of~$S$
and thus closed.

(b) is a special case of\,(a).
\end{proof}
Having verified the closed product property, we obtain:
%
%\ma{JH2}
\begin{thm}\label{JH2}
Let $G$ be a totally disconnected,
locally compact group,\linebreak
$\alpha\in \Aut(G)$ be contractive,
and $\Omega:=\langle \alpha\rangle$
or $\Omega:=\langle\Int(G)\cup \{\alpha\}\rangle$.
Then $G$ admits an
$\Omega$-composition series,
and the
Schreier Refinement Theorem
and the Jordan-H\"{o}lder Theorem
hold in the form described in Proposition~{\rm \ref{preJordH}}.
\end{thm}
\begin{proof}
The existence of an $\Omega$-composition series
follows from Lemma~\ref{easybit} below.
The
remainder holds by Proposition~\ref{preJordH}
and Corollary~\ref{havecpp}\,(b).
\end{proof}
We complete the proof using a well-known
fact (cf.\ \cite[Proposition III.13.20]{FaD}):
%
%\ma{fromHaR}
\begin{la}\label{fromHaR}
Let $\alpha$ be an automorphism of
a locally compact group $G$
and $N$ be an $\alpha$-stable closed
normal subgroup of~$G$.
Let $\wb{\alpha}$ be the automorphism
induced by $\alpha$ on $G/N=:Q$.
Then
$\Delta_G(\alpha)=\Delta_N(\alpha|_N)
\cdot\Delta_Q(\wb{\alpha})$.\Punkt
\end{la}
%
%
%\ma{easybit}
\begin{la}\label{easybit}
Let $(G,\alpha)$ be a totally disconnected contraction group.
Then the length of any $\langle\alpha\rangle$-series
$\one=G_0\triangleleft \,\cdots\, \triangleleft G_n=G$
without repetitions
is bounded by the number of prime factors
of $\Delta_G(\alpha^{-1})$, counted with
multiplicities.
\end{la}
\begin{proof}
Let $\alpha_j\colon G_j/G_{j-1}\to G_j/G_{j-1}=:Q_j$
be the contractive automorphism of~$Q_j$
induced by~$\alpha$.
Then
$\Delta_G(\alpha^{-1})\!=\!
\Delta_{Q_1}(\alpha_1^{-1})\Delta_{Q_2}(\alpha_2^{-1})\cdots
\Delta_{Q_n}(\alpha_n^{-1})$
by Lemma~\ref{fromHaR},
where $1<\Delta_G(\alpha^{-1}_j)\in \N$
for each $j\in \{1,\ldots, n\}$,
by Proposition~\ref{basefacts}\,(e).
The assertion is now immediate.
\end{proof}
For later use,
we record another important consequence of Proposition~\ref{inflate}.
%
%\ma{pivotal}
\begin{cor}\label{pivotal}
Let $(G,\alpha)$ and $(H,\beta)$
be totally disconnected, locally compact
contraction groups
and $\phi\colon G\to H$ be a continuous
homomorphism such that $\beta\circ\phi=\phi \circ\alpha$.
Then $\phi(G)$ is
$\beta$-stable,
closed in~$H$,
and $\phi|^{\phi(G)}\colon
G\to\phi(G)$ is a quotient map.
In particular,
if $\phi$ is injective, then~$\phi$ is
a topological isomorphism onto its image.
\end{cor}
\begin{proof}
Let $W\leq G$ be a compact open subgroup
such that $\alpha(W)\sub W$.
Because $\phi(W)$ is a compact subgroup
of~$H$ with $\beta(\phi(W))=\phi(\alpha(W))\sub \phi(W)$,
Proposition~\ref{inflate}\,(a) shows that
$S:=\bigcup_{n\in \N_0}\beta^{-n}(\phi(W))$
is a closed subgroup of~$H$
which possesses $\phi(W)$ as a compact, open subgroup.
Since $G=\bigcup_{n\in \N_0}\alpha^{-n}(W)$
and $\phi\circ\alpha=\beta\circ\phi$,
we deduce that $\phi(G)=S$.
Because $W\sub G$ is compact,
$\beta|_W^{\beta(W)}$
is a quotient morphism of topological
groups and hence open.
Since $\phi(W)$ is open in~$S$,
we deduce that also $\phi|^S \colon G\to S$
is an open map, which completes the proof.
\end{proof}
\begin{rem}
Using Corollary~\ref{pivotal},
all standard facts concerning
Remak decompositions of finite groups,
as formulated in \cite[3.3.1--3.3.10]{Rob},
can be adapted directly to totally
disconnected contraction groups,
notably the Krull-Remak-Schmidt Theorem.
The corollary ensures that all homomorphisms
encountered in the classical proofs are continuous,
and all isomorphisms bicontinuous.
\end{rem}
\section{Simple contraction groups are pro-discrete}\label{secprodis}
In this section, we show that every non-trivial
contraction group $(G,\alpha)$
has a non-trivial $\alpha$-stable closed
normal subgroup $S\triangleleft G$ which is pro-discrete.
Therefore every simple contraction group
is pro-discrete.
This information is essential
for the proof of the classification.
%
%
%
%\ma{givesprodis}
\begin{prop}\label{givesprodis}
Let $(G,\alpha)$ be a totally disconnected
contraction group.
\begin{itemize}
\item[\rm (a)]
Then $G$ has a largest
closed
normal $\alpha$-stable
subgroup~$S^\alpha(G):=S$
possessing
a compact, open, $\alpha$-invariant
subgroup which is normal in~$G$.
If $G$ is non-trivial, then
also $S$ is non-trivial.
\item[\rm (b)]
$S$ can be obtained as follows:
Let $W\leq G$ be a compact, open subgroup
such that $\alpha(W)\sub W$,
and
$N:=\bigcap_{x\in G}xWx^{-1}$
be its core.
Then $N$ is an $\alpha$-invariant
closed normal subgroup of~$G$
and $S=\bigcup_{n\in \N_0}\alpha^{-n}(N)$.
\end{itemize}
\end{prop}
{\bf Proof.}
We may assume without loss of generality
that $G\not=\one$.
Let $W\leq G$ and $N$ be as in (b);
then clearly $N$ is closed,
and it is the largest normal
subgroup of~$G$ contained in~$W$.
Furthermore,
$\alpha(N)=
\bigcap_{x\in G}x\alpha(W)x^{-1}\sub N$.
Thus $S:=\bigcup_{n\in \N_0}\alpha^{-n}(N)$
is a closed
normal subgroup of~$G$
possessing~$N$ as an open subgroup,
by Proposition~\ref{inflate}\,(a)
and\,(d).
If $H\triangleleft G$ is any
$\alpha$-stable closed normal
subgroup of~$G$
possessing an $\alpha$-invariant,
compact, open subgroup $K\leq H$ which is
normal in~$G$,
then there is $m\in \N$ such that
$\alpha^m(K)\sub W$
and thus $\alpha^m(K)\sub N$
(since $\alpha^m(K)\triangleleft G$),
entailing that
$H=\bigcup_{n\in \N_0}\alpha^{-n}(\alpha^m(K))\sub
\bigcup_{n\in \N_0}\alpha^{-n}(N)=S$.
Thus, it only remains to prove
that~$S\not=\one$. We proceed in steps.

%\ma{helpnormalize}
\begin{numba}\label{helpnormalize}
For each $n\in \N_0$, the set $V_n:=\bigcap_{x\in W}x\alpha^n(W)x^{-1}$
is compact, and it is the largest
subgroup of~$\alpha^n(W)$ which is normal
in~$W$.
Since $O:=\alpha^n(W)$ is open in~$W$
and normalizes $\alpha^n(W)$,
we see that $V_n=\bigcap_{xO\in W/O}x\alpha^n(W)x^{-1}$
is open in~$W$.
For each $k\in \N_0$,
we have
%
%\ma{givesinsight}
\begin{equation}\label{givesinsight}
\alpha^{-k}(V_n)\;=\;
\bigcap_{x\in \alpha^{-k}(W)}x\alpha^{n-k}(W)x^{-1}\,.
\end{equation}
Since $W\sub \alpha^{-k}(W)$,
we see that $\alpha^{-k}(V_n)$
is normalized by~$W$, for each $k\in \N_0$.
\end{numba}
\begin{numba}
\emph{$\alpha^{-k}(V_n)$ is
a normal subgroup of~$W$, for all
$n\in \N_0$ and $k\in \{0,1,\ldots, n\}$.}
Indeed,
we have $\alpha^{n-k}(W)\sub W$
and thus $\alpha^{-k}(V_n)\sub W$,
by~(\ref{givesinsight}).
As $W$ normalizes $\alpha^{-k}(V_n)$,
the assertion follows.
\end{numba}
%
%
%\ma{preparr}
\begin{numba}\label{preparr}
\emph{$\alpha^{-1}(V_n)\not\sub V_n$ holds,
for each $n\in \N_0$.}
Otherwise
$\alpha^{-k}(V_n)\sub V_n$
and thus
%
%
%\ma{givecontra}
\begin{equation}\label{givecontra}
V_n\;\sub \; \alpha^k(V_n)\,, \quad
\mbox{for all $k\in \N_0$.}
\end{equation}
Because $V_n\not=\one$, there exists an identity
neighbourhood $P\sub G$ which is a proper
subset of~$V_n$. Then $\alpha^k(V_n)\sub P$
for large~$k$, since $\alpha$ is
compactly contractive.
This contradicts (\ref{givecontra}).
\end{numba}
%
%
%\ma{combto}
\begin{numba}\label{combto}
\emph{For each $n\in \N_0$, we have
$U_n:=\alpha^{-n}(V_n)\not\sub \alpha(W)$.}
Indeed, otherwise
$\alpha^{-1}(V_n)=\alpha^{n-1}(U_n)\sub
\alpha^n(W)\sub W$.
Hence $\alpha^{-1}(V_n)$
is a subgroup of $\alpha^n(W)$ which is
normal in~$W$ (by {\bf\ref{helpnormalize}}).
As $V_n$ is the largest such subgroup
(see {\bf \ref{helpnormalize}}),
we have $\alpha^{-1}(V_n)\sub V_n$.
This contradicts~{\bf\ref{preparr}}.
\end{numba}
\begin{numba}
Since $U_n=\bigcap_{x\in \alpha^{-n}(W)}xWx^{-1}$,
where $W\sub \alpha^{-1}(W)\sub \alpha^{-2}(W)\sub \cdots$
and $G=\bigcup_{n\in \N_0}\alpha^{-n}(W)$,
we see that $U_1\supseteq U_2\supseteq \cdots$ and
%
%\ma{UUn}
\begin{equation}\label{UUn}
\bigcap_{n\in \N_0}U_n\;=\:
\bigcap_{x\in G}xWx^{-1}\;=\; N\,.
\end{equation}
Since $U_n\cap (W\setminus \alpha(W))\not=\emptyset$
by {\bf\ref{combto}},
the set
$\{U_n\cap (W\setminus \alpha(W))\colon
n \in \N_0\}$
of compact sets has the finite
intersection property, and thus
$N\cap (W\setminus \alpha(W))=
\bigcap_{n\in \N_0}(U_n\cap(W\setminus \alpha(W)))
\not=\emptyset$,
showing that~$N\not=\one$
and hence also $S\not=\one$.
This completes the proof.\,\,\Punkt
\end{numba}
Note that $S^\alpha(G)$ is pro-discrete,
by Proposition~\ref{basefacts}\,(f).
We readily deduce:
%
%\ma{simplarepro}
\begin{cor}\label{simplarepro}
Every simple, totally disconnected
contraction group is pro-discrete.\Punkt
\end{cor}
\section{Further technical tools}\label{sectools}
In this section, we compile
two technical lemmas,
which will be used to prove the classification.
The first of these provides information
concerning the closed, normal,
$\alpha$-invariant subgroups
of a simple contraction group.
%
%\ma{convenien}
\begin{la}\label{convenien}
Let $(G,\alpha)$ be a simple totally disconnected
contraction group.
If $N\triangleleft G$ is
a closed normal subgroup
such that $\alpha(N)\sub N$,
then either $N=\one$ or $N$ is open in~$G$.
\end{la}
\begin{proof}
By Proposition~\ref{inflate}\,(a) and\,(d),
$S:=\bigcup_{k\in \N_0}\alpha^{-k}(N)$
is an $\alpha$-stable, closed normal
subgroup of~$G$ possessing~$N$
as an open subgroup.
Since $(G,\alpha)$ is simple
we either have $S=\one$ or $S=G$.
The assertion follows.
\end{proof}
The next lemma will
be used later to identify those
simple contraction groups which are abelian
torsion groups.
Here and in the following,
two contraction groups
$(G,\alpha)$ and $(H,\beta)$
are called \emph{isomorphic}
if there exists an isomorphism of
topological groups $\phi\colon G\to H$
such that $\beta\circ\phi=\phi\circ\alpha$.
\begin{la}\label{createshift}
Let $(G,\alpha)$ be a simple totally disconnected
contraction group.
If there exists a non-trivial,
finite, normal subgroup~$N\triangleleft G$,
then $(G,\alpha)$ is isomorphic
to $F^{(-\N)}\times F^{\N_0}$
with the right shift,
for some finite, simple group~$F$.
\end{la}
\begin{proof}
Let $F\sub N$ be a minimal non-trivial
normal subgroup of~$G$.
For $n\in \N_0$, consider the map
\[
\phi_n\colon F^{\{0,1,\ldots, n\}} \to G\,,\quad (x_0,\ldots, x_n)\mto
x_0\alpha(x_1)\cdots \alpha^n(x_n)\,.
\]
We show by induction on $n\in \N_0$
that $\phi_n$ is an injective
homomorphism.
This is trivial if $n=0$.
If $\phi_n$ is an injective homomorphism,
then $F\alpha(F)\cdots \alpha^n(F)$
is a normal subgroup of~$G$,
whence also its image
$N:=\alpha(F)\cdots \alpha^{n+1}(F)$
is a normal subgroup of~$G$.
Hence either $F\cap N=\one$
(entailing that indeed
$\phi_{n+1}$ is an injective homomorphism),
or $F\cap N=F$ and thus $F\sub N$,
by minimality of~$F$.
If $F\sub N$, then
$\alpha^{-1}(N)=
F\alpha(F)\cdots \alpha^n(F)
\sub N\alpha(F)\cdots \alpha^n(F)=N$
and thus $\alpha(N)=N$,
as $N$ is finite.
Hence $\alpha^k(N)=N$ for each
$k\in \N$, contradicting the
fact that $\alpha^k(N)\to\one$ as
$\alpha$ is compactly contractive.

There is a compact, open, normal subgroup $W\triangleleft G$
such that $\alpha(W)\sub W$,
and a maximal number $k\in \Z$ such that
$F\sub \alpha^k(W)$.
Then
$F\cap \alpha^{k+1}(W)$
is a normal subgroup of~$G$ and a proper subset of~$F$
(by maximality of~$k$), and thus
$F\cap \alpha^{k+1}(W)=\one$
by minimality of~$F$.
Hence, after replacing $W$ with $\alpha^k(W)$,
without loss of generality
$F\sub W$ and $F\cap \alpha(W)=\one$.
We define
\[
\phi\colon F^{\N_0}\to G\,,\quad
\phi((x_n)_{n\in \N_0})\, :=\,
\lim_{n\to\infty} \phi_n(x_0,\ldots, x_n)\,;
\]
the limits exist because
$\phi_n(x_0,\ldots,x_n)^{-1}\phi_{n+m}(x_0,\ldots,
x_{n+m})\in \alpha^{n+1}(W)$
for all $n,m\in \N_0$.
Each $\phi_n$ being a homomorphism,
also~$\phi$ is a homomorphism.
Given $m\in \N_0$,
the set $U_m:=\{(x_n)_{n\in \N_0}\in F^{\N_0}\colon
\mbox{$x_n=1$ for all $n<m$}\}$
is an identity neighbourhood
in $F^{\N_0}$, and
$\phi(U_m)\sub \alpha^m(W)$.
Therefore $\phi$ is continuous at~$1$
and hence continuous.
Furthermore, $\phi$ is injective.
To see this, let
$x=(x_n)_{n\in \N_0}\in F^{\N_0}$
such that $x\not=1$.
There exists a smallest
integer $m\in \N_0$ such that
$x_m\not=1$.
Then $\phi(x)=\alpha^m(x_m)y$
with $y=\lim_{n\to\infty}\alpha^{m+1}(x_{m+1})\cdots
\alpha^{m+n}(x_{m+n})\in \alpha^{m+1}(W)$.
If $\phi(x)=1$, then
$x_m^{-1}=\alpha^{-m}(y)\in F\cap \alpha(W)=\one$
and thus $x_m=1$, which is a contradiction.
Thus $\phi(x)\not=1$, whence $\ker\phi=\one$
and~$\phi$ is injective.
The image~$K$ of~$\phi$ is compact,
and it is normal in~$G$, being the closure
of the normal subgroup $\bigcup_{n\in \N_0}\im\phi_n$.
Furthermore, $\alpha(K)\sub K$
and $K\not=\one$.
Hence~$K$ is open in~$G$, by Lemma~\ref{convenien}.
Set $H:=F^{(-\N)}\times F^{\N_0}$,
$V:=F^{\N_0}\sub H$
and let $\sigma\colon H\to H$ be the
right shift.
Then $\sigma(V)\sub V$
and $\phi\circ \sigma|_V=\alpha\circ \phi$,
by construction of~$\phi$.
As $\alpha$ and $\sigma$ are contractive
automorphisms,
\cite[Proposition~2.2]{Wan} shows that
$\phi$ extends to
a topological isomorphism $H\to G$
that intertwines $\sigma$ and~$\alpha$.
\end{proof}
%
%
%
%
%
%
%
%
%
%
%
%\ma{seccompfac}
\section{Proof of the classification}\label{seccompfac}
We are now in the position to
prove the classification of the simple
totally disconnected contraction groups
described in Theorem~A.\vspace{.5mm}
\begin{center}
{\bf Classification of the abelian simple contraction
groups}\vspace{-.5mm}
\end{center}
We first determine
a system of representatives for
the abelian simple contraction groups.
For a discussion of automorphisms
of totally disconnected, locally compact
abelian groups, see also~\cite{ERW}.
%
%
%\ma{theabcase}
\begin{thm}\label{theabcase}
Let $(G,\alpha)$ be a simple
totally disconnected contraction group
which is abelian.
Then $G$ is locally
pro-$p$ for some prime~$p$
and $(G,\alpha)$
is either isomorphic to $C_p^{(-\N)}\times C_p^{\N_0}$
with the right shift,
or isomorphic to
$(\Q_p^d,\beta)$ for some $d\in \N$
and a contractive $\Q_p$-linear automorphism
$\beta\colon \Q_p^d\to\Q_p^d$
which does not leave any proper non-trivial
vector subspace invariant.
\end{thm}
\begin{proof}
Let $W\leq G$ be a compact, open subgroup
such that $\alpha(W)\sub W$.
Then $W$ has a non-trivial
$p$-Sylow subgroup $W_p$
for some prime~$p$.
Since~$W$ is abelian, $W_p$
is unique
(cf.\ \cite[Proposition 2.2.2\,(d)]{Wsn}).
Since $\alpha(W_p)\leq W$ is
a pro-$p$-group, we have $\alpha(W_p)\sub W_p$
(cf.\ \cite[Proposition 2.2.2\,(b)]{Wsn}).
By Lemma~\ref{convenien}, $W_p$ is open
in~$G$. Hence $G$ is locally pro-$p$,
and we may assume that $W=W_p$.
Then $[W:\alpha(W)]=p^m$
for some $m\in \N$.
For each $n\in \N$, we
write $W^{\{n\}}:=\{x^n\colon
x\in W\}$.
Since $W^{\{p\}}$
is a closed, normal,
$\alpha$-invariant subgroup of~$G$,
either $W^{\{p\}}=\one$ or $W^{\{p\}}$ is open
in~$G$, by Lemma~\ref{convenien}.
If $W^{\{p\}}=\one$,
then $W$ (and hence $G$)
is a torsion group of exponent~$p$,
whence $(G,\alpha)$ is isomorphic to
$C_p^{(-\N)}\times C_p^{\N_0}$,
as a consequence of Lemma~\ref{createshift}.
If $W^{\{p\}}$
is an open subgroup of~$G$,
then $\alpha^n(W)\sub W^{\{p\}}$
for some $n\in \N$.
Thus
\begin{equation}
W^{\{p^k\}}\; \supseteq \; \alpha^{nk}(W)\quad
\mbox{for all $k\in \N$,}
\end{equation}
by induction:
we have $W^{\{p^{k+1}\}}=(W^{\{p^k\}})^{\{p\}}
\supseteq (\alpha^{nk}(W))^{\{p\}}=\alpha^{nk}(W^{\{p\}})
\supseteq \alpha^{nk}(\alpha^n(W))=\alpha^{n(k+1)}(W)$.
Thus
$[W:W^{\{p^k\}}]\leq[W:\alpha^{nk}(W)]=$\linebreak
$[W:\alpha(W)]\cdots [\alpha^{nk-1}(W):\alpha^{nk}(W)]
=[W:\alpha(W)]^{nk}=p^{nmk}$
and so
%
%\ma{givesfinrk}
\begin{equation}\label{givesfinrk}
[W:W^{\{p^k\}}]\leq p^{nmk}\quad\mbox{for each $k\in \N$.}
\end{equation}
Using \cite[Theorem\,3.16]{DMS},
we deduce from
(\ref{givesfinrk})
that the pro-$p$-group
$W$ has finite rank.
Therefore $G$ is a $p$-adic
Lie group (see \cite[Corollary~8.33]{DMS}).
Since $G$ is an abelian Lie group,
its Lie algebra $L(G)$
is an abelian Lie algebra.
Therefore the Campbell-Hausdorff
multiplication coincides
with the addition map $L(G)\times L(G)\to L(G)$,
and we find an exponential
map $\phi\colon P\to Q$
which is an isomorphism
of topological groups from
a compact, open additive subgroup
$P\leq L(G)$ onto a compact, open subgroup
$Q\leq G$.
Set $\beta:=L(\alpha)\colon L(G)\to L(G)$.
After shrinking~$P$ and~$Q$,
we may assume that
%
%\ma{naturalty}
\begin{equation}\label{naturalty}
\phi\circ \beta|_U=\alpha\circ \phi|_U
\end{equation}
for some open subgroup $U\sub P$ such
that $\beta(U)\sub P$.
After shrinking $U$, we may assume
that $V:=\phi(U)=\alpha^N(W)$
for some $N\in \N$.
Since $\alpha(V)\sub V$,
we deduce from (\ref{naturalty})
that $\beta(U)\sub U$.
Hence $\phi\circ \beta^n|_U=\alpha^n\circ \phi|_U$
for each $n\in \N$,
entailing that $\beta$ is a contractive automorphism
of $L(G)$.
Now $(G,\alpha)$ is equivalent
to $(L(G),\beta)$ by \cite[Proposition~2.2]{Wan}.
\end{proof}
It remains to describe normal forms
in the $p$-adic case.
\begin{defn}
Given a prime~$p$,
let $R_p \sub \Q_p[X]$ be the set of all
irreducible monic polynomials~$f$ whose
roots in an algebraic closure
$\wb{\Q_p}$ of~$\Q_p$
have absolute value $<1$.
For $f=X^d+a_{d-1}X^{d-1}+\cdots+a_0 \in R_p$,
we set $E_f:=\Q^d$,
let $e_1,\ldots, e_d\in \Q^d$
be the standard basis vectors
and define $\alpha_f$ as the linear
automorphism of~$E_f$
determined by $\alpha_f(e_j)=e_{j+1}$
for $j\in \{1,\ldots, d-1\}$
and $\alpha_f(e_d)=-\sum_{i=1}^d a_{i-1} e_i$.
\end{defn}
Note that $R_p$ has continuum cardinality,
as $\{X\!-\! a\colon a\in \Q_p, |a|<1\}\sub R_p$.
%
%\ma{ratform}
\begin{prop}\label{ratform}
The family
$(E_f,\alpha_f)_{f\in R_p}$
is a system of representatives
for the isomorphism classes
of the simple totally disconnected contraction
groups $(G,\alpha)$
such that $G$ is abelian,
torsion-free, and locally pro-$p$.
\end{prop}
\begin{proof}
Abbreviate $\K:=\Q_p$
and let $\wb{\K}$ be an algebraic closure
of~$\Q_p$.
Given $(G,\alpha)$ as described in the proposition,
we may assume that $G=E$
is a finite-dimensional $\K$-vector space
and $\alpha$ a continuous linear map,
by Proposition~\ref{theabcase}.
We consider $E$ as a $\K[X]$-module
via $X.v:=\alpha(v)$ for $v\in E$.
Then $E$ is irreducible
and thus $E\isom \K[X]/f\K[X]$
for a unique monic
irreducible polynomial
$f\in \K[X]$
(cf.\ \cite[\S3.9, Exercise~2]{Jc1}).
Given $r>0$, let
$\wb{E}_r$ be the sum of all generalized
eigenspaces of $\alpha\tensor \wb{\K}$
in $E\tensor_\K\wb{\K}$
to eigenvalues $\lambda\in \wb{\K}$
such that $|\lambda|=r$.
Then $E=\bigoplus_{r>0}E_r$,
where $E_r:=\wb{E}_r\cap E$
is $\alpha$-stable
(see \cite[p.\,81]{Mar}),
and thus $E=E_r$ for some $r>0$,
as $(E,\alpha)$ was assumed simple.
There exists an ultrametric norm $\|.\|$
on $E=E_r$ such that $\|\alpha(x)\|=r\|x\|$
for each $x\in E$ (see \cite{SPO}).
Since $\alpha$ is contractive,
it follows that $r<1$
and thus $f \in R_p$.
Let $d$ be the degree of $f$.
With respect to the basis corresponding to
$X^0,\ldots, X^{d-1}$,
the automorphism $\alpha$ has the
same matrix as $\alpha_f$
with respect to $e_1,\ldots, e_d$,
and thus $(G,\alpha)$ is isomorphic
to $(E_f,\alpha_f)$.

Conversely, for each $f\in R_p$,
the $\K[X]$-module $E_f$ (with $X.v:=\alpha_f(v)$)
is irreducible
and uniquely determines~$f$
(cf.\ \cite[\S3.9]{Jc1}).
By irreducibility,
$E_f=(E_f)_r$ for some $r>0$
(with notation as before),
where $r<1$ by definition of~$R_p$.
As there exists an ultrametric norm~$\|.\|$
on $E_f=(E_f)_r$ such that $\|\alpha_f(v)\|=r\|v\|$
for each $v\in E_f$,
we see that $\alpha_f$ is a contractive automorphism.
To complete the proof,
let $N\sub E_f$ be a non-trivial,
$\alpha_f$-stable closed additive subgroup.
Then $\Spann_{\Q_p}(N)$
is an $\alpha_f$-stable, non-trivial
vector subspace of~$E_f$ and hence
$E_f=\Spann_{\Q_p}(N)$, by irreducibility.
Since $N$ is open in $\Spann_{\Q_p}(N)=E_f$
and $\alpha_f$ is contractive,
we deduce that
$E_f=\bigcup_{n\in \N_0}\alpha^{-n}_f(N)=N$.
Thus $(E_f,\alpha_f)$
is a simple contraction group.
\end{proof}
By the preceding proof,
for each $f\in R_p$
all eigenvalues $\lambda$
of $\alpha_f$ in $\wb{\Q_p}$ have the same absolute
value $r:=|\lambda|$.\vspace{.5mm}
\begin{center}
{\bf Classification of the non-abelian simple contraction
groups}\vspace{-.5mm}
\end{center}
To classify the non-abelian simple
contraction groups, we shall use a folklore lemma
from group theory (the proof of which
is recalled in Appendix~\ref{appdetails}).
%
%
%\ma{abstrgp}
\begin{la}\label{abstrgp}
Let $G$ be a group and $N_1,\ldots, N_n$
be pairwise distinct normal subgroups of~$G$
such that $G/N_k$ is a non-abelian simple
group, for each\linebreak
$k\in \{1,\ldots, n\}$.
Abbreviate $D:=N_1\cap\cdots\cap N_n$
and $D_k:=\bigcap_{j\not=k}N_j$.
Then
\[
\theta\colon G/D\to G/N_1\times\cdots\times G/N_n\,,
\quad xD\mto (xN_1,\ldots, xN_n)
\]
is an isomorphism of groups
which takes $D_k/D$ isomorphically onto
$G/N_k$, for each $k\in \{1,\ldots, n\}$.\Punkt
\end{la}
%
%
%
%\ma{nonabcase}
\begin{thm}\label{nonabcase}
Let $(G,\alpha)$
be a simple totally disconnected contraction group.
If $G$ is non-abelian,
then $(G,\alpha)$ is isomorphic to
$F^{(-\N)}\times F^{\N_0}$
with the right shift
for a non-abelian, finite simple group~$F$.
\end{thm}
\begin{proof}
We let $W\triangleleft G$ be a compact, open, normal
subgroup such that $\alpha(W)\sub W$.
As $G$ is non-abelian,
there are $g,h\in G$ such that
$ghg^{-1}h^{-1}\not=1$.
After applying a suitable power of $\alpha$ two both
elements, we may assume that $h\in W$.
There is $m \in \N$ such that
$ghg^{-1}h^{-1}\not\in \alpha^m(W)$.
As a consequence,
$g\not\in \ker\phi$
for the homomorphism
\[
\phi\colon G\to \Aut(W/\alpha^m(W))\,,\quad
\phi(x)(w\alpha^m(W))\, :=\,
xwx^{-1}\alpha^m(W)\,.
\]
Since $\alpha^m(W)\sub \ker\phi$,
we deduce that $N:=\ker\phi$ is a proper,
open, normal subgroup of~$G$.
The group $\Aut(W/\alpha^m(W))$ being finite,
$N$ has finite index in~$G$.
Consequently,
there exists a maximal proper normal
subgroup $M$ of~$G$
such that $N\sub M$.
Then $M$ is open, and $F:=G/M$
is a finite simple group.
The closure $\wb{[G,G]}$
of the commutator subgroup is a
non-trivial, $\alpha$-stable
closed normal subgroup of~$G$
and thus $\wb{[G,G]}=G$,
by simplicity.
Hence $[F,F]=F$ and thus~$F$ is non-abelian.
Next, we observe that
%
%\ma{resfin}
\begin{equation}\label{resfin}
\bigcap_{k\in \Z}\alpha^k(M)\;=\; \one\,,
\end{equation}
because this intersection is an $\alpha$-stable,
closed, normal, proper subgroup of~$G$
and $(G,\alpha)$ is simple.
Here $\alpha^k(M)$ is
a normal subgroup of~$G$
such that $G/\alpha^k(M)\isom G/M=F$
is a non-abelian, simple
group.
If $k_1\not=k_2$, say $k_2>k_1$,
then $\alpha^{k_1}(M)\not=\alpha^{k_2}(M)$
because otherwise
$\alpha^{k_2-k_1}(M)=M$,
entailing that $\bigcap_{k\in \Z}\alpha^k(M)
=\bigcap_{k=0}^{k_2-k_1-1}\alpha^k(M)$
is open as a finite intersection
of open sets, which is absurd.
Let $\psi\colon G\to G/M=F$
be the quotient homomorphism
and set $\psi_k:=\psi\circ \alpha^{-k}$
for $k\in \Z$.
Since $\ker \psi_k=\alpha^k(M)$,
in view of the properties
just established Lemma~\ref{abstrgp}
shows that the map
\[
\psi_{n,m}\, :=\,
(\psi_k)_{k=n}^m\colon G\to F^{\{n,\ldots, m\}}\,,\quad
x\mto (\psi_n(x),\ldots, \psi_m(x))
\]
is surjective,
for all $n,m\in \Z$ such that $n\leq m$.
Given $x\in G$, there is $k_0\in \Z$ such that
$x\in \alpha^k(M)$ for all $k\leq k_0$.
We can therefore define
a homomorphism
\[
\eta:=(\psi_k)_{k\in \Z}
\colon G\to F^{(-\N)}\times F^{\N_0}
\,,\quad
\eta(x)\, :=\, (\psi_k(x))_{k\in \Z}\,.
\]
We let
$\sigma$ be the right shift on
$F^{(-\N)}\times F^{\N_0}=:H$.
Then $\sigma\circ \eta=\eta\circ \alpha$
by construction of~$\eta$.
To complete the proof,
we show that~$\eta$ is
an isomorphism of topological groups.
First,
$\eta$ is injective,
because $\ker\eta=\one$ by (\ref{resfin}).
Since $G=\bigcup_{n\in \Z}\alpha^{-n}(W)$
as an ascending union,
there is $n\in \Z$ such that
$\alpha^{-n}(W)\not\sub M$.
On the other hand,
$\alpha^k(W)\sub M$
for large~$k$
since~$\alpha$
is compactly
contractive.
Hence $n$ can be chosen
minimal.
As a consequence,
$W\sub\ker\psi_k$
for all $k<n$
while $W\not\sub \ker\psi_k$
for all $k\geq n$.
Thus $\eta(W)\sub F^{\{n,n+1,\ldots\}}$.
Since~$H$
induces the product topology
on $F^{\{n,n+1,\ldots\}}$,
we see that $\theta:=\eta|_W$
is continuous and hence also~$\eta$.
Let $m\geq n$. 
Because $W$ is normal in~$G$
and $\psi_{n,m}\colon G\to F^{\{n,\ldots, m\}}$
is surjective, the image $\psi_{n,m}(W)$
is a normal subgroup of
the product $F^{\{n,\ldots, m\}}$
of non-abelian simple groups.
By Remak's Theorem \cite[3.3.12]{Rob},
$\psi_{n,m}(W)=F^J$
for a subset $J\sub \{n,\ldots, m\}$.
Since $\psi_k(W)\not={\bf 1}$
for each $k\in \{n,\ldots,m\}$,
we see that $J=\{n,\ldots, m\}$
and thus $\psi_{n,m}(W)=F^{\{n,\ldots, m\}}$.
As a consequence, $\theta$ has dense
image. Now $\theta(W)$ being also
compact and thus closed, we deduce that
$\eta(W)=\theta(W)=F^{\{n,n+1,\ldots\}}$.
Hence $\eta$ has open image,
and since $\psi_{n,m}$ is
surjective for all $n,m\in\Z$ such that
$n\leq m$, we see that $\eta(G)$
is dense in~$H$ and hence equal to~$H$.
Because $\eta|_W$ is a homeomorphism
onto its open image, $\eta$ is an isomorphism
of topological groups.
\end{proof}
\begin{rem}
The finite group $F$ in 
Theorem~\ref{nonabcase} is uniquely determined
up to isomorphism.
To see this, note that every
compact, open, normal subgroup
$W\sub F^{(-\N)}\times F^{\N_0}$
such that $\sigma(W)\sub W$
is of the form $W=F^{\{n,n+1,\ldots\}}$
for some $n\in \Z$,
and $W/\sigma(W)\isom F$.
\end{rem}
\section{Canonical {\boldmath $\alpha$}-stable series}\label{sec:series}
We now describe how a series
$\one =S^\alpha_0(G) \triangleleft
S^\alpha_1(G)\triangleleft\cdots \triangleleft S^\alpha_n(G) =G$
of $\alpha$-stable closed normal subgroups
can be associated
to each totally disconnected contraction
group $(G,\alpha)$ in a canonical way.
This series will serve as a technical tool
in the proof of the Structure Theorem
(Theorem~B).
It can also be used
to see that $\langle \Int(G)\cup\{\alpha\}\rangle$-composition
factors are pro-discrete
(Proposition~\ref{alsoprodis}).
\begin{defn}
Let $G$ be a totally disconnected, locally compact
group and $\alpha\colon G\to G$ be a contractive
automorphism.
We define $S^\alpha_0(G):=\one$
and
$S^\alpha_1(G):=S^\alpha(G)$ (as in
Proposition~\ref{givesprodis}\,(a)).
Inductively,
having defined the $\alpha$-stable
closed normal subgroup
$S_{j-1}^\alpha(G)\triangleleft G$,
we set $Q_j:=G/S^\alpha_{j-1}(G)$,
let $q_j\colon G\to Q_j$ be the
quotient map and set
$S_j^\alpha(G):=q^{-1}_j(S^{\alpha_j}(Q_j))$,
where $\alpha_j \colon Q_j\to Q_j$
is the contractive automorphism
determined by $\alpha_j\circ q_j=q_j\circ \alpha$.
By Lemma~\ref{easybit},
the series $S^\alpha_0(G)\triangleleft
S^\alpha_1(G)\triangleleft \cdots$ becomes stationary.
Thus, there is a smallest
$n\in \N_0$
such that $S^\alpha_{n+1}(G)=S^\alpha_n(G)$.
We call
$\one=S^\alpha_0(G)\triangleleft S^\alpha_1(G)\triangleleft \cdots\triangleleft
S_n^\alpha(G)=G$ the \emph{canonical
$\alpha$-stable series of~$G$}.
\end{defn}
Let us call
an $\langle\Int(G)\cup\{\alpha\}\rangle$-series
$\one=G_0\triangleleft G_1\triangleleft \cdots\triangleleft G_n=G$
\emph{special}
if $G_i/G_{i-1}$ has a compact, open
subgroup which is normal in $G/G_{i-1}$
and invariant under the
automorphism of $G/G_{i-1}$
induced by~$\alpha$,
for each $i\in \{1,\ldots, n\}$.
The following proposition compiles
various useful properties
of the canonical $\alpha$-stable series.
In particular, it is special
and ascends faster than any other special
$\langle\Int(G)\cup\{\alpha\}\rangle$-series.
%
%
%\ma{propmore}
\begin{prop}\label{propmore}
Let $(G,\alpha)$ be a totally disconnected contraction
group and $\Omega:=\langle \Int(G)\cup\{\alpha\}\rangle$.
Then the following holds:
\begin{itemize}
\item[\rm (a)]
The canonical
$\alpha$-stable series
$\one =S^\alpha_0(G) \triangleleft
S^\alpha_1(G)\triangleleft\cdots \triangleleft S^\alpha_n(G) =G$
is a special $\Omega$-series
without repetitions.
\item[\rm (b)]
$S^\alpha_j(G)/S^\alpha_{j-1}(G)$
is a pro-discrete, closed
normal subgroup of $Q_j:=G/S^\alpha_{j-1}(G)$,
for each $j\in \{1,\ldots, n\}$.
\item[\rm (c)]
If $\one=G_0\triangleleft G_1\triangleleft \cdots \triangleleft  G_m=G$
is any special
$\Omega$-series,
then $m\geq n$ and $G_j\sub S^\alpha_j(G)$
for each $j\in \{0,\ldots, n\}$.
\item[\rm (d)]
If $\beta\in \Aut(G)$
such that $\beta\circ\alpha=\alpha\circ\beta$,
then $\beta(S^\alpha_j(G))=S^\alpha_j(G)$
for all $j\in \{0,\ldots, n\}$.
\end{itemize}
\end{prop}
\begin{proof}
(a) and (b):
By construction,
each $S^\alpha_j(G)$ is an $\alpha$-stable, closed
normal subgroup of~$G$
and hence an
$\Omega$-subgroup.
If $G\not=\one$, then
$S^\alpha_1(G)\not=\one$ by
Proposition~\ref{givesprodis}\,(a)
and thus $\one=S^\alpha_0(G)\subset
S^\alpha_1(G)$.
Likewise
$S^{\alpha_j}(Q_j)\not={\bf 1}$ for $j\in \{1,\ldots, n\}$
and hence
$S_{j-1}^\alpha(G)\subset q_j^{-1}(S^{\alpha_j}(Q_j))=
S^\alpha_j(G)$, using the notation of
the preceding definition. Therefore
no repetitions occur.
By Proposition~\ref{givesprodis}\,(a),
$S_j^\alpha(G)/S_{j-1}^\alpha(G)=S^{\alpha_j}(Q_j)$
has an $\alpha_j$-invariant,
compact, open subgroup~$W_j$
which is normal in~$Q_j$.
Hence the canonical $\alpha$-stable series
is a special $\Omega$-series.
In particular, $W_j$
is a compact, open,
normal subgroup of
$S_j^\alpha(G)/S_{j-1}^\alpha(G)=S^{\alpha_j}(Q_j)$,
and hence
$S_j^\alpha(G)/S_{j-1}^\alpha(G)$
is pro-discrete by
Proposition~\ref{basefacts}\,(f).

(c) We show by induction on $i\in \{0,\ldots, n\}$
that $m\geq i$ and $G_i\sub S^\alpha_i(G)$.
For $i=0$, this is clear.
Now assume that $i\in \{0,\ldots, n\}$,
$m\geq i-1$, and $G_{i-1}\sub S^\alpha_{i-1}(G)$.
Since $S^\alpha_{i-1}(G)$ is a proper subset of~$G$,
by the preceding so is $G_{i-1}$
and hence $m\geq i$.
By hypothesis,
there exists a compact, open subgroup
$K\sub G_i/G_{i-1}$ which is normal
in $G/G_{i-1}$
and invariant under the automorphism of $G/G_{i-1}$
induced by~$\alpha$.
The continuous homomorphism $q\colon
G/G_{i-1}\to G/S_{i-1}^\alpha(G)$,
$xG_{i-1}\mto xS_{i-1}^\alpha(G)$
intertwines $\alpha$
and the contractive automorphism
$\alpha'$
of $G/S_{i-1}^\alpha(G)$
induced by~$\alpha$.
As a consequence of Corollary~\ref{pivotal},
$q(G_i/G_{i-1})$ is a closed,
$\alpha'$-stable normal subgroup
of $G/S_{i-1}^\alpha(G)$
which has $q(K)$ as an open
subgroup.
Since $q(K)$ is normal in $G/S_{i-1}^\alpha(G)$
and $\alpha'$-invariant,
Proposition~\ref{givesprodis}\,(a)
shows that $q(G_i/G_{i-1})
\sub S^{\alpha'}(G/
S_{i-1}^\alpha(G))
=S_i^\alpha(G)/S_{i-1}^\alpha(G)$.
Therefore $G_iS_{i-1}^\alpha(G)\sub S_i^\alpha(G)$
and thus $G_i\sub S_i^\alpha(G)$.

(d) If $\beta\in \Aut(G)$ commutes
with $\alpha$, then $\beta(S^\alpha_1(G))$
is an $\alpha$-stable closed
normal subgroup of~$G$ containing an $\alpha$-invariant,
compact open subgroup which is normal in~$G$,
and thus $\beta(S^\alpha_1(G))\sub S^\alpha_1(G)$.
Likewise, $\beta^{-1}(S^\alpha_1(G))\sub S^\alpha_1(G)$,
and thus $\beta(S^\alpha_1(G))=S^\alpha_1(G)$.
The same argument can be applied to
$Q_j=G/S^\alpha_{j-1}(G)$ and its automorphisms
$\alpha_j$ and $\beta_j$ induced
by $\alpha$ and $\beta$;
hence a simple induction yields
the assertion.
\end{proof}
We know from
Corollary~\ref{simplarepro}
that composition factors of $\langle\alpha\rangle$-composition
series are pro-discrete.
As a first application of the canonical
$\alpha$-stable series, we now show that also
$\langle \Int(G)\cup\{\alpha\}\rangle$-composition
factors are pro-discrete.
%
%
%
%
%\ma{alsoprodis}
\begin{prop}\label{alsoprodis}
Let $(G,\alpha)$ be a totally disconnected
contraction group and $\Omega:=\langle \Int(G)\cup\{\alpha\}\rangle$.
Then the factors of each
$\Omega$-composition series
of $(G,\alpha)$ are pro-discrete.
Furthermore,
$(G,\alpha)$ has
an $\Omega$-composition series
which is a special $\Omega$-series.
\end{prop}
\begin{proof}
Since all $\Omega$-composition series
are equivalent by the Jordan-H\"{o}lder Theorem
(Theorem~\ref{JH2}), to prove the
first assertion
it suffices to consider an $\Omega$-composition series
$\one = H_0 \triangleleft H_1\triangleleft \cdots \triangleleft H_m=G$
which has been obtained by
refining the canonical
$\alpha$-stable series of~$G$
(this is possible
by the Schreier Refinement Theorem).
Let $i\in \{1,\ldots, m\}$.
Then $i\in \{k+1,\ldots, k+\ell\}$
for some $k\in \{0,\ldots, m-1\}$
and $\ell\in \{1,\ldots, m-k\}$
such that, for some $j\in \{1,\ldots, n\}$,
$H_k=S_{j-1}^\alpha(G)$ and
$H_{k+\ell}=S_j^\alpha(G)$.
By Proposition~\ref{propmore}\,(a),
$H_{k+\ell}/H_k= S_j^\alpha(G)/S_{j-1}^\alpha(G)$
has a compact, open subgroup~$W/H_k$
which is normal in $G/H_k$
and invariant under the automorphism
of $G/H_k$ induced by~$\alpha$.
Now standard arguments show that
$(W\cap H_i)/H_{i-1}$
is a compact, open subgroup
of $H_i/H_{i-1}$ which is normal
in $G/H_{i-1}$
and invariant under the automorphism
of $H_i/H_{i-1}$ induced by~$\alpha$.
\end{proof}
\section{Proof of the Structure Theorem}\label{sec:outline}
We now outline
the main steps of the proof of the Structure Theorem\linebreak
(Theorem~B from the Introduction).
The details of Steps 2 to 4
will be given
in Sections~\ref{secextprod} to
\ref{sec:divisible}.\\[2.5mm]
Throughout the following, $(G,\alpha)$
is a totally disconnected, locally compact
contraction group (unless we state the contrary).
Furthermore,
%
%
%\ma{currentserie}
\begin{equation}\label{currentserie}
{\mathbf 1} \,= \,G_0 \; \triangleleft \; \cdots\;
\triangleleft \; G_n \, = \, G
\end{equation}
is an $\langle\alpha\rangle$-composition
series for $G$.
By the classification,
each factor $G_j/G_{j-1}$ is pro-discrete and is isomorphic to
either
(a) $(\Q_p^d,+)$ for some prime $p$ and some $d\in \N$; or to
(b) a restricted product over $\Z$ of copies
of a finite simple group.
In case~(a), $G_j/G_{j-1}$ is infinitely divisible
and torsion-free.
In case~(b), $G_j/G_{j-1}$ is a torsion group
of finite exponent.\\[2.5mm]
It is useful to consider the special cases
first where either all $\langle\alpha\rangle$-composition factors
are torsion groups, or all of them are torsion-free.\\[5mm]
{\bf Step 1. The case when all composition factors
are torsion groups.}\\[1.5mm]
If each of the factors $G_j/G_{j-1}$
is a torsion group of finite exponent,
then also $G$ is a torsion group
of finite exponent,
as a special case of the following lemma
(the proof of which is based
on obvious inductive arguments):
%
%
%\ma{obviou}
\begin{la}\label{obviou}
Let $\one=G_0\triangleleft
G_1\triangleleft \cdots\triangleleft G_n=G$
be a series of groups.
\begin{itemize}
\item[\rm (a)]
If $G_j/G_{j-1}$ is a torsion
group for each $j\in \{1,\ldots, n\}$,
then $G$ is a torsion group.
If $G_j/G_{j-1}$ is a torsion
group of exponent $m_j$
for each $j\in \{1,\ldots, n\}$,
then $G$ is a torsion group
of finite exponent which divides
$m_1\cdot\ldots\cdot m_n$.
\item[\rm (b)]
If $G_j/G_{j-1}$ is torsion-free
for $j=k, \ldots, n$,
then $\tor(G)=\tor(G_{k-1})$.\,\,\Punkt\vspace{3mm}
\end{itemize}
\end{la}
The next lemma implies
that the exponent of a torsion factor
$G_j/G_{j-1}$\linebreak
divides the module of the automorphism
induced by $\alpha^{-1}$ on $G_j/G_{j-1}$.
This information will be useful later.
%
%\ma{scaleshift}
\begin{la}\label{scaleshift}
Let $F$ be a finite group,
$H:=F^{(-\N)}\times
F^{\N_0}$
and $\sigma$ be the right shift on~$H$.
Then
$\Delta_H(\sigma^{-1})=|F|$.
\end{la}
\begin{proof}
For the compact open subgroup
$W:=F^{\N}$ of~$H$, we have
$\sigma^{-1}(W)=F^{\N_0}$
and $\Delta_H(\sigma^{-1})=\lambda(\sigma^{-1}(W))/\lambda(W)
=[\alpha^{-1}(W):W]=|F|$.
\end{proof}
{\bf Step~2. The special case of torsion-free composition factors.}\\[1.5mm]
The following proposition
(proved in Section~\ref{secextprod})
describes the structure
of contraction groups all of whose composition factors
are torsion-free.
%
%
%\ma{struc-t-free}
\begin{prop}\label{struc-t-free}
Let $(G,\alpha)$ be a totally disconnected
contraction group possessing an
$\langle\alpha\rangle$-composition series
$\one=G_0\triangleleft
%G_1\triangleleft
\cdots\triangleleft G_n=G$
such that $G_j/G_{j-1}$ is torsion-free,
for each $j\in \{1,\ldots, n\}$.
Then $G$ is an internal direct product
$G=G_{p_1}\times\cdots\times G_{p_r}$
of certain nilpotent $p$-adic Lie groups~$G_p$.
Each $G_p$ is topologically fully invariant
in~$G$ $($and hence $\alpha$-stable$)$.
\end{prop}
In the situation of Proposition~\ref{struc-t-free},
$G$ is infinitely divisible and torsion-free,
as a consequence of the next lemma.
%
%
%\ma{div-torfree}
%
\begin{la}\label{div-torfree}
Let $G$ be a $p$-adic Lie group admitting a contractive
automorphism~$\alpha$.
Then $G$ is infinitely divisible and torsion-free.
\end{la}
\begin{proof}
Let $\exp\colon V\to U$ be an exponential map
of~$G$, which is a diffeomorphism
from an open $\Z_p$-submodule $V\sub L(G)$
onto an open subgroup $U\leq G$.
Then each $x\not=1$ in $U$ has the form $x = \exp(X)$ for some
$X\not=0$ in~$V$. For each $n\in\N$, we then have
$x^n = \exp(nX) \ne 1$, showing
that $U$ is torsion-free
and hence also $G=\bigcup_{k\in \N_0}\alpha^{-k}(U)$.
Furthermore, $\{x^n\colon x\in U\}=\exp(nV)$
is an identity neighbourhood in~$G$
consisting of elements possessing an $n$-th root.
Hence every element of $G=\bigcup_{k\in \N_0}
\alpha^{-k}(\exp(nV))$ has an $n$-th root.\vspace{2mm}
\end{proof}
{\bf Step~3. The set of torsion elements
is a subgroup.}\\[1.5mm]
If one of the composition factors in (\ref{currentserie})
is torsion,
then it may be assumed that $G$ has torsion elements
and that $G_1$ is torsion (see Section~\ref{sec:torsion}).
As a consequence, it may always be supposed that
torsion factors appear first
in the composition series:
%
%\ma{re-org}
%
\begin{la}\label{re-org}
Each totally disconnected contraction group $(G,\alpha)$
admits an $\langle\alpha\rangle$-composition series
${\mathbf 1} = G_0 \triangleleft \cdots
\triangleleft G_n = G$ such that, for suitable\linebreak
$k\in \{0,\ldots, n\}$,
the factors $G_j/G_{j-1}$ are torsion groups
for $j\in\{1,\ldots, k\}$
and all other factors are torsion-free.
Then $\tor(G)=G_k$ is a subgroup of~$G$.
\end{la}
The proof of Lemma~\ref{re-org}
uses the following result. 
%
%
%\ma{torchar}
%
\begin{la}\label{torchar}
If $\one = G_0 \triangleleft \cdots \triangleleft G_k$
is an $\langle\alpha\rangle$-composition series for
$G_k$ with $G_i/G_{i-1}$
a torsion group for $i\in \{1,\ldots, j\}$ and
$G_i/G_{i-1}$ a torsion-free group for $i\in \{j+1,\ldots, k\}$,
then $G_j=\tor(G_k)$ is a characteristic
subgroup of $G_k$ and $G_k/G_j$ is torsion-free.
\end{la}
\begin{proof}
Lemma~\ref{obviou}\,(a) and (b) show that
$G_j$ is a torsion group and $\tor(G_k)=\tor(G_j)=G_j$.
Thus $G_j$ is a characteristic subgroup of~$G_k$.
By Proposition~\ref{struc-t-free}
and Lemma~\ref{div-torfree},
$G_k/G_j$
is torsion-free.
\end{proof}
Now Lemma~\ref{re-org} readily follows:
If $n=0$ or if all factors $G_j/G_{j-1}$
are torsion-free, or if all factors
are torsion, there is nothing to show.
Now assume that
$n$ is arbitrary and that $G$ has torsion
elements but is not a torsion group.
By Lemma~\ref{lemmins},
we may assume that $G_1/G_0$ is torsion,
whence there exist $k,m\in \{1,\ldots, n\}$
with $m>k$
such that $G_j/G_{j-1}$ is torsion
for all $j\in \{1,\ldots, k\}$
while $G_j/G_{j-1}$ is torsion-free
for $j\in \{k+1,\ldots, m\}$
and~$m$ cannot be increased.
We assume that the $\langle\alpha\rangle$-composition
series has been chosen such that~$k$ is maximal.
Then $G_k=\tor(G_m)$ by Lemma~\ref{obviou}.
If $m=n$, there is nothing more to show.
Otherwise, $G_{m+1}/G_m$ is a torsion group
and since $G_k=\tor(G_m)$ is characteristic
in $G_m$, we deduce that $G_k$ is normal in $G_{m+1}$.
Now the torsion factor in the composition series
of $G_{m+1}/G_k$ can be moved to
the bottom, and hence $G_{k+1}/G_k$ can be replaced by a torsion
factor, contradicting the maximality of~$k$.\vspace{4.5mm}\Punkt

\noindent
{\bf Step~4. Definition of a complementary subgroup {\boldmath $D$}.}\\[1.5mm]
We now choose the $\langle\alpha\rangle$-composition series
(\ref{currentserie}) as described in Lemma~\ref{re-org}.
Thus $G_j/G_{j-1}$ is torsion
for $j\in \{1,\ldots, k\}$
while $G_j/G_{j-1}$ is torsion-free
and infinitely
divisible for $j\in \{k+1,\ldots, n\}$.
Then $T:=\tor(G)=G_k$
is a characteristic (and hence $\alpha$-stable)
subgroup of~$G$.
%
%
%\ma{power}
%
\begin{la}\label{power}
Put $t_\alpha := \Delta_T(\alpha^{-1}|_T)$.
Then  $x^{t_\alpha} = 1$ for all $x\in T$.
\end{la}
\begin{proof}
Consider the
$\langle \alpha\rangle$-series
$\one =G_0\triangleleft G_1 \triangleleft \cdots \triangleleft G_k=T$.
For $j\in\{1,\dots,k\}$,
let $\alpha_j$ be the automorphism induced by $\alpha$ on
$Q_j:=G_j/G_{j-1}$
and put $t_j := \Delta_{Q_j}(\alpha_j^{-1})$.
Then $t_\alpha = t_1\cdots t_k$
(see proof of Lemma~\ref{easybit}).
Furthermore, $Q_j$ is a torsion group
of exponent dividing~$t_j$ (cf.\ Lemma~\ref{scaleshift}).
Thus $T$ is a torsion group
of finite exponent that divides~$t_\alpha$,
by Lemma~\ref{obviou}\,(a).
\end{proof}
Define $D := \wb{\langle x^{t_\alpha} \mid x\in G \rangle}$.
Then $D$ is a closed, topologically
characteristic (and hence $\alpha$-stable)
subgroup of $G$.
We record
an essential property of~$D$:
%
%
%\ma{surjective}
%
\begin{la}\label{surjective}
The map $\phi \colon D\to G/T$, $\phi(x):=xT$
is surjective.
\end{la}
\begin{proof}
Let $yT$ be in $G/T$, where $y\in G$.
Since, by Lemma~\ref{div-torfree}, $G/T$ is infinitely
divisible, there is $xT\in G/T$ such that
$yT = (xT)^{t_\alpha}$.
Then $yT = x^{t_\alpha}T$ belongs to the range of $\phi$.
\end{proof}
The following lemma (established
in Section~\ref{sec:divisible})
completes the proof of the first half of
the Structure Theorem:
%
%\ma{factor}
\begin{la}\label{factor}
$D$ is an infinitely divisible group,
$D=\dv(G)$, and $G = T\times D$.
\end{la}
Now also the second half of the Structure Theorem
readily follows:
Since $D\cap T=\{1\}$,
the group~$D$ is torsion-free.
Hence all composition factors
of $D$ are torsion-free (see Section~\ref{sec:torsion})
and therefore $D=G_{p_1}\times\cdots\times G_{p_r}$
is an internal direct product
of $\alpha$-stable $p$-adic Lie groups~$G_p$,
by Proposition~\ref{struc-t-free}.
%
%
%
%
%
%
%
%
%
%
%\ma{secextprod}
%
\section{The case of torsion-free factors}\label{secextprod}
In this section, we prove Proposition~\ref{struc-t-free},
thus completing Step~2 of Section~\ref{sec:outline}.
The proof is based on the following lemma.
%
%
%
%\ma{info-prod}
\begin{la}\label{info-prod}
Let $N$ be a closed normal subgroup
of a topological group~$G$
and $Q:=G/N$.
\begin{itemize}
\item[\rm (a)]
If $N$ and $Q$ are $p$-adic Lie groups,
then $G$ is a $p$-adic Lie group.
\item[\rm (b)]
If $Q$ is a $q$-adic Lie group
and $N$ an internal direct product
$N=\prod_{p\in \cp}N_p$ of $p$-adic
Lie groups $N_p$,
where $\cp$ is a finite set of primes,
then $G$ has an open subgroup~$U$ that is an internal
direct product $U=\prod_{p\in \cp\, \cup\{q\}}U_p$
of $p$-adic Lie groups~$U_p$.
\item[\rm (c)]
If $N$ is $\alpha$-stable
for a contractive automorphism~$\alpha$ of~$G$
in the situation of {\rm (b)},
then $G$ is an internal direct product
$G=\prod_{p\in \cp \cup\{q\}}G_p$
of $\alpha$-stable $p$-adic Lie groups~$G_p$.
\end{itemize}
\end{la}
\begin{proof}
(a) $G$ is locally compact
by \cite[(5.25)]{HaR}
and totally disconnected,
whence it has a compact
open subgroup~$U$.
Then $N\cap U$ and $U/(N\cap U)\isom UN/N\sub G/N$
are $p$-adic Lie groups,
and after shrinking~$U$ both of these
groups are pro-$p$-groups
of finite rank (by \cite[Corollary~8.33]{DMS}
and \cite[Proposition~8.1.1\,(a)]{Wsn}).
By \cite[Proposition~1.11\,(ii)]{DMS}
and \cite[Proposition~8.1.1\,(b)]{Wsn},
also $U$ is a pro-$p$-group
of finite rank
and thus $G$ is a $p$-adic Lie group by~\cite[Corollary~8.33]{DMS}.

(b) We may assume that $q\in \cp$
(otherwise, define $N_q:=\one$).
For each $p\in \cp$, we let $V_p\sub N_p$
be an open subgroup which is a pro-$p$-group
(see \cite[Corollary~8.33]{DMS}).
Let $U\sub G$ be a compact, open subgroup
such that $U\cap N\sub \prod_{p\in \cp}V_p$
and such that $\pi(U)$ is a pro-$q$-group,
where $\pi\colon G\to Q$ is the quotient map.
It readily follows from
\cite[Proposition~2.4.3]{Wsn}
that $U\cap N=\prod_{p\in \cp}(U\cap V_p)$
(see \cite[Proposition~2.2]{MIX});
hence $U\cap N=\prod_{p\in \cp}V_p$
without loss of generality.
Being the unique $p$-Sylow subgroup,
$V_p$ is topologically characteristic in the normal
subgroup $U\cap N$ of~$U$,
and hence $V_p$ is a normal subgroup of~$U$.
Given $p\in \cp\setminus\{q\}$,
let $U_p$ be a $p$-Sylow subgroup of~$U$.
Then $\pi(U_p)=\one$ and thus
$U_p\sub U\cap N$,
entailing that $U_p=V_p$ is a $p$-adic
Lie group and a normal subgroup of~$U$.
Hence $M:=\prod_{p\not=q}U_p$ is
a normal subgroup of~$U$.
Let $U_q$ be a $q$-Sylow subgroup of~$U$
containing $V_q$ (see \cite[Proposition~2.2.2\,(c)]{Wsn}).
Then $M\cap U_q=\one$
and $MU_q$ is a subgroup
of~$U$.
By \cite[Proposition~2.2.3\,(b)]{Wsn}, $\pi(MU_q)=\pi(U_q)$
is a $q$-Sylow subgroup of~$\pi(U)$
and thus $\pi(MU_q)=\pi(U)$.
As $MU_q$ is saturated
under $U\cap N=MV_q$,
we deduce that $U=MU_q$.
Thus $U=M\semid \, U_q$.
For $p\not=q$,
the conjugation action of
$U_q$ on~$U_p$ gives rise
to a continuous homomorphism
$\phi_p\colon U_q\to\Aut(U_p)$,
where $\Aut(U_p)$
is equipped with the compact-open topology
(cf.\ \cite[Theorem~3.4.1]{Eng}).
Like every compact $p$-adic Lie group,
$U_p$ is topologically finitely generated.
Hence $\Aut(U_p)$ has an open
subgroup $O_p$ which is a pro-$p$-group
(see \cite[Theorem~5.6]{DMS}
and \cite[Theorem~4.4.2]{RaZ}).
Since every continuous homomorphism from
a pro-$q$-group to a pro-$p$-group is trivial,
we deduce that $\phi_p^{-1}(O_p)\sub\ker\phi_p$.
Hence $W_q:=\bigcap_{p\not=q}\phi_p^{-1}(O_p)$
is an open subgroup of~$U_q$.
After replacing $U$
with its open subgroup $MW_q$
and $U_q$ with $W_q$,
we may assume that $U_q$
centralizes $U_p$ for each
$p\not=q$. Thus $U=M\times U_q=\prod_{p\in \cp}U_p$
as an internal direct product.
It only remains to observe that $U_q$
is a $q$-adic Lie group by~(a),
because $\pi(U_q)=\pi(U)$ and $N\cap U_q=V_q$
are $q$-adic Lie groups.

(c) Without loss of generality $q\in \cp$.
By (b), $G$ has a compact open subgroup
$U$ which is a direct product $U=\prod_{p\in \cp}U_p$
of $p$-adic
Lie groups.
After shrinking $U$, we may assume that
each $U_p$ is a pro-$p$-group
and hence topologically fully invariant
in~$U$ (being its unique $p$-Sylow subgroup).
There exists a compact, open subgroup
$W\leq G$ such that $W\sub U$ and
$\alpha(W)\sub W$.
Then $W=\prod_{p\in\cp}(W\cap U_p)$;
after replacing $U$ with $W$,
we may assume that $\alpha(U)\sub U$
and thus $\alpha(U_p)\sub U_p$
for each $p\in \cp$.
By Proposition~\ref{inflate}\,(a),
$G_p:=\bigcup_{k\in \N_0}\alpha^{-k}(U_p)$
is an $\alpha$-stable
closed subgroup of~$G$
which has $U_p$ as an open subgroup.
Hence $G_p$ is a $p$-adic Lie group.
Let $p_1,\ldots,p_r$ be the distinct elements of~$\cp$
and consider the product map
$\psi\colon \prod_{p\in \cp}G_p\to G$,
$(x_{p_1},\ldots,x_{p_r})\mto x_{p_1}\cdots x_{p_r}$.
Then $\beta:=\alpha|_{G_{p_1}}\times\cdots\times
\alpha|_{G_{p_r}}$ is a contractive automorphism of
$\prod_{p\in \cp}G_p$
and~$\psi$ intertwines $\beta$ and $\alpha$,
\begin{equation}\label{intertw}
\alpha\circ \psi \; = \; \psi\circ\beta\, .
\end{equation}
We now use that $\psi$ induces a bijection
from $\prod_{p\in \cp}U_p$ onto~$U$.
Since $\psi$ is an injective homomorphism
on $\prod_{p\in \cp}U_p$
and $\bigcup_{k\in \N_0}\beta^{-k}(\prod_{p\in \cp}U_p)=\prod_{p\in \cp}G_p$,
we deduce from (\ref{intertw})
that $\psi$ is injective
and a homomorphism.
Hence $G_p\triangleleft G$ for each~$p$.
Since $U\sub\im(\psi)$
and $\bigcup_{k\in \N_0}\alpha^{-k}(U)=G$,
using (\ref{intertw})
we see that $\psi$ is also surjective. Hence
$\psi$ is an isomorphism.
\end{proof}
{\bf Proof of Proposition~\ref{struc-t-free}.}
Each composition factor $G_j/G_{j-1}$ being a $p$-adic
Lie group for some prime~$p$
by Theorem~A, a straightforward
induction on~$n$ based on
Part\,(c) of Lemma~\ref{info-prod}
shows that $G$ is an internal direct product
$G=G_{p_1}\times\cdots\times G_{p_r}$
of certain $\alpha$-stable
subgroups~$G_p$
which are $p$-adic Lie groups.
Each $G_p$ has an open pro-$p$
subgroup $U_p$;
then $\alpha^{-k}(U_p)$
also is a pro-$p$-group
for each $k\in \N$, and
$G_p=\bigcup_{k\in \N_0}\alpha^{-k}(U_p)$.
Since, for $p\not=q$,
each continuous homomorphism
from a pro-$p$-group
to a pro-$q$-group is trivial,
we deduce that each
endomorphism of the topological
group $G$ takes $G_p$ to $G_p$.
Thus $G_p$ is topologically fully invariant.
To complete the proof,
we recall
from \cite[Theorem~3.5\,(ii)]{Wan}
that every $p$-adic contraction group
is a unipotent $p$-adic algebraic group
and hence nilpotent.\Punkt
%
%
%
%
%
%
%
%
%\ma{sec:torsion}
\section{Shifting torsion factors to the bottom}\label{sec:torsion}
In this section, we prove the following
lemma, which completes the details
of Step~3 in Section~\ref{sec:outline}.
%
%\ma{lemmins}
\begin{la}\label{lemmins}
Let $(G,\alpha)$ be a totally
disconnected contraction group
such that at least one $\langle\alpha\rangle$-composition factor
of~$G$ is a torsion group.
Then $G$ has an $\langle\alpha\rangle$-composition series
{\rm (\ref{currentserie})}
such that $G_1$ is a torsion group.
In particular, $G$ has non-trivial
torsion elements.
\end{la}
{\bf Proof.}
If the lemma was false,
we could find a counterexample
with an $\langle\alpha\rangle$-composition
series
of minimal length $n\geq 2$.
By minimality,
for each $\langle \alpha\rangle$-composition series
$\one = G_0
%\triangleleft G_1
\triangleleft \cdots \triangleleft G_n= G$,
the factors $G_j/G_{j-1}$ have to be
torsion-free for all $j\in \{1,\ldots, n-1\}$,
while $G_n/G_{n-1}$ is a torsion group.
By Proposition~\ref{struc-t-free},
$G_{n-1}=H_{p_1}\times\cdots\times H_{p_r}$
is a direct product of certain topologically
characteristic (and hence $\alpha$-stable)
nilpotent $p$-adic Lie groups~$H_p\not=\{1\}$.
If $r\geq 2$,
then $K:=H_{p_2}\times\cdots\times H_{p_r}$
is topologically characteristic in $G_{n-1}$
and hence normal in~$G_n$.
We may assume that $K=G_j$ for some $j\in \{1,\ldots, n-2\}$.
Because $Q:=G/K$ has a properly shorter
$\langle\alpha\rangle$-composition series
than $G$,
it has an $\langle\alpha\rangle$-composition
series starting in a torsion factor
and thus also $G_{j+1}/G_j$
can be chosen as a torsion group,
which is a contradiction.
Thus $G_{n-1}$ is a nilpotent $p$-adic
Lie group for some~$p$.
The closed commutator subgroup
$C$ of $G_{n-1}$ being topologically
characteristic in $G_{n-1}$,
arguing as before we reach a
contradiction unless $C={\bf 1}$.
Hence $G_{n-1}$ is an abelian $p$-adic
Lie group. The next two lemmas
will help us to reach a final contradiction.
%
%
%\ma{inserted-la}
\begin{la}\label{inserted-la}
Let $(G,\alpha)$ be a totally disconnected
contraction group and \mbox{$N\triangleleft G$} be an
$\alpha$-stable closed normal subgroup.
If $G/N$ is a torsion group
and $N$ an abelian $p$-adic Lie group,
then $N$ is contained in the centre of~$G$.
\end{la}
\begin{proof}
As in the proof of Proposition~\ref{theabcase},
we see that
$N\isom \Q_p^d$ for some $d\in \N$.
Hence $\Aut(N)\isom \GL_d(\Q_p)$
is a $p$-adic Lie group (when equipped with
the compact-open topology).
Therefore $\Aut(N)$
has a torsion-free open subgroup~$W$.
Since $N$ is abelian, a
homomorphism of groups can be defined~via
\[
\phi\colon G/N\to\Aut(N)\,,\qquad
\phi(xN)(y)\, :=\, xyx^{-1}\,.
\]
To see that $\phi$ is continuous,
let $U\leq G$ be a compact, open
subgroup; then $V:= UN/N$ is a compact, open subgroup
of $Q:=G/N$
and we only need to show that
$\phi$ is continuous on~$V$.
By \cite[Proposition~1.3.3]{Wsn},
there exists a continuous
section $\sigma\colon V\to U$
to the quotient morphism
$U\to V$, $x\mto xN$
of pro-finite groups.
Then the map
\[
V\times N\to N\,,\qquad
(x,y)\mto\phi(x)(y)=\sigma(x)y\sigma(x)^{-1}
\]
is continuous and hence so is
$\phi|_V\colon V \to\Aut(N)$
(cf.\ \cite[Theorem~3.4.1]{Eng}).
Now $\phi$ being continuous,
$\phi^{-1}(W)$ is an identity neighbourhood
in~$Q$. Since $Q$ is a torsion group
and $W$ is torsion-free,
we must have $\phi^{-1}(W)\sub\ker \phi$
and thus $\ker\phi$ is open.
Let $\bar{\alpha}$ be the contractive automorphism
of $Q$ induced by~$\alpha$.
Given $xN\in \ker\phi$,
we have
$\phi(\bar{\alpha}^{-1}(xN))(\alpha^{-1}(y))
=\alpha^{-1}(x)\alpha^{-1}(y)\alpha^{-1}(x^{-1})
=\alpha^{-1}(\phi(xN)(y))=\alpha^{-1}(y)$
for each $y\in N$,
and thus $\bar{\alpha}^{-1}(xN)\in\ker\phi$.
Hence $Q=\bigcup_{k\in \N_0}\bar{\alpha}^{-k}(\ker\phi)=\ker\phi$.
As a consequence,
$xyx^{-1}=\phi(x N)(y)=y$ for each $x\in G$
and thus $N\sub Z(G)$.
\end{proof}
%
%
%
%\ma{enterjacob}
\begin{la}\label{enterjacob}
Let $(G,\alpha)$ be a totally disconnected
contraction group, $A\sub G$
be a central, $\alpha$-stable closed subgroup
and $q\colon G\to Q$ be a quotient morphism
with kernel~$A$.
Assume that
$A\isom \Q_p^d$
and $Q\isom F^{(-\N)}\times F^{\N_0}$
with the right shift~$\sigma$,
for a finite group~$F$.
Then
$\tor(G)$ is a subgroup of~$G$,
and $G=A\times \tor(G)$ internally as a topological group.
\end{la}
\begin{proof}
Without loss of generality $Q=F^{(-\N)}\times F^{\N_0}$.
We set $F_k:=F^{\{{-k},\ldots, k\}}$ for $k\in \N$,
$G_k:=q^{-1}(F_k)$,
and consider $A$ as an $F_k$-module with
the trivial action.
For each $n\in \N$,
the $n$-th cohomology group
$H^n(F_k,A)$ with coefficients in~$A$
(as in \cite[\S6.9]{Jac})
is a $\Q_p$-vector space in a natural
way and hence a torsion-free group.
On the other hand, $F_k$ being finite,
$H^n(F_k,A)$ is a torsion group by \cite[Theorem~6.14]{Jac}.
Hence $H^n(F_k,A)=\{0\}$ for each
$n\in\N$ and thus $H^2(F_k,A)=\{0\}$
in particular, entailing that the extension
$A\to G_k\to F_k$ splits
and thus $G_k=A\times S_k$ internally
for some subgroup $S_k\leq G_k$
(cf.\ \cite[Theorem~6.15]{Jac}).
Since $A$ is torsion-free and
$S_k\isom F_k$ a torsion group,
we deduce that $S_k=\tor(G_k)$.
In particular, $S_k$ is uniquely determined
and $S_k=\tor(G_k)\sub \tor(G_{k+1})=S_{k+1}$,
whence $S_\infty:=\bigcup_{k\in \N}S_k$
and its closure $S:=\wb{S_\infty}$
are subgroups of~$G$.
Each $S_k\isom F_k$ being a torsion
group of exponent dividing~$|F|$,
also $S_\infty$ and its closure $S$
are torsion groups
of exponent dividing~$|F|$ and thus $S\cap A=\one$,
because $A$ is torsion-free.

From
$q(\alpha(G_k))=\sigma(q(G_k))=\sigma(F_k)\sub F_{k+1}$
we deduce that $\alpha(G_k)\sub G_{k+1}$
and thus $\alpha(S_k)=\alpha(\tor(G_k))\sub \tor(G_{k+1})=S_{k+1}$,
entailing that $\alpha(S_\infty)\sub S_\infty$.
Likewise, $\alpha^{-1}(S_\infty)\sub S_\infty$,
whence $\alpha(S_\infty)=S_\infty$
and $\alpha(S)=S$.
By Corollary~\ref{pivotal}, $q(S)$ is closed in~$Q$.
Since $q(S)\supseteq q(S_\infty)
=\bigcup_{n\in \N}F_n=F^{(\Z)}$,
where $F^{(\Z)}$ is dense in~$Q$, we see that
$q(S)=Q$. Hence $G=A\times S$ internally as a group
and hence also as a topological group
(cf.\ Corollary~\ref{pivotal}). Since $A$ is torsion-free
and $S\isom Q$ a torsion group,
$S=\tor(G)$ follows.
\end{proof}
\emph{Proof of Lemma~{\rm \ref{lemmins}}, completed.}
Since $G_n/G_{n-1}$ is a torsion group,
and is a factor
of an $\langle\alpha\rangle$-composition series,
it is isomorphic to a restricted
product of copies of a finite group.
Lemmas~\ref{inserted-la} and \ref{enterjacob} show that $G_n$ has a
closed subgroup isomorphic to this restricted product.
This subgroup is characteristic in~$G_n$
and can be chosen as $G_1$. We have reached a
contradiction.\Punkt
\section{Proof that {\boldmath $D$}
has the desired properties}\label{sec:divisible}
In this section, we prove Lemma~\ref{factor},
thus completing the details of Step~4
from Section~\ref{sec:outline}.\\[2.5mm]
To prove Lemma~\ref{factor},
we use induction on the length of the canonical
$\alpha$-stable series of~$G$;
the induction starts because the case
$G=\one$ is trivial. For general~$G$, let
$S:= S^\alpha_1(G)$ be the first term in the canonical
series for~$G$.
Then $G/S$ has a shorter canonical series and so,
by the inductive hypothesis, $G/S = \tilde T \times \tilde D$,
where $\tilde T$ is the torsion subgroup of $G/S$
and $\tilde D
:=\wb{\langle x^{\tilde{t}_\alpha}\colon x\in G/S\rangle}$
is an infinitely divisible subgroup,
where $\tilde{t}_\alpha$ is the module
of the automorphism of $\tilde T$ induced by $\alpha^{-1}$.
Let $q \colon G \to G/S$ be the quotient map.
Then $T\sub q^{-1}(\tilde T)$. 

By Proposition~\ref{givesprodis}\,(a), there is a compact,
open subgroup $N\leq S$
such that $N\triangleleft G$, $\alpha(N) \sub N$
and $S = \bigcup_{j\in {\mathbb Z}} \alpha^{-j}(N)$.
Since $\alpha$ is an automorphism, $\alpha^j(N)\triangleleft G$
for every $j\in{\mathbb N}$.
Consider, for each $j\in {\mathbb N}$,
the finite group $N/\alpha^j(N)$ and define a
homomorphism 
\[
\phi_j \colon
G \to \hbox{Aut}(N/\alpha^j(N))\
\hbox{ by }\  \phi_j(x) \colon
y\alpha^j(N) \mapsto xyx^{-1}\alpha^j(N) \ \ \mbox{for $y\in N$.}
\]
Then $\hbox{ker}(\phi_j)$ is a finite index normal subgroup
of $G$ for each $j$ and so there is a positive integer,
$d_j$, such that $x^{d_j} \in \hbox{ker}(\phi_j)$ for
every $x\in G$. Since $\alpha$ is an automorphism, we have
that $x^{d_j} \in \alpha^\ell(\hbox{ker}(\phi_j))$ for every
$x\in G$ and $\ell\in {\mathbb Z}$.
Define $M_j :=
\bigcap_{\ell\in{\mathbb Z}} \alpha^\ell(\hbox{ker}(\phi_j))$.
Then $M_j$ is a
closed $\alpha$-stable normal subgroup
of $G$ for each $j$. 

It is clear that
$\hbox{ker}(\phi_{j+1})\sub \hbox{ker}(\phi_j)$ for each~$j$.
Hence $(M_j)_{j\in{\mathbb N}}$ is a decreasing sequence
of closed $\alpha$-stable normal subgroups of~$G$.
As $\Delta_{M_{j+1}}(\alpha^{-1}|_{M_{j+1}})$ is
a positive integer strictly less than $\Delta_{M_j}
(\alpha^{-1}|_{M_j})$ if $M_{j+1}$ is a proper normal
subgroup of $M_j$ (by Proposition~\ref{basefacts}\,(e)
and Lemma~\ref{fromHaR}),
this sequence eventually stabilizes.
Thus there is a $J$ such that $M_j = M_J$ for all $j\geq J$. 
%
%
%\ma{commute}
\begin{la}\label{commute}
$M_J$ is equal to the centralizer of $S$.
\end{la}
\begin{proof}
It is clear from the definition of $M_j$ that the
centralizer of $S$ is a subgroup of $M_j$ for every $j$.

For the converse, let $s\in S$. By the definition of~$N$,
there is an
$\ell\in{\mathbb Z}$ such that $s\in \alpha^\ell(N)$.
If $x\in M_J$, then for every $j\geq J$ we have
$xsx^{-1}\alpha^{\ell+j}(N) = s\alpha^{\ell+j}(N)$.
Since $\bigcap_{j\geq J} \alpha^{\ell+j}(N) = \one$
(because $\alpha$ is compactly contractive),
it follows that $xsx^{-1} = s$. 
\end{proof}
The subgroup $M_J$ is not trivial if $G$ has a
torsion-free composition factor. 
%
%
%\ma{MJnottrivial}
\begin{la}\label{MJnottrivial}
The map $\phi \colon D\cap M_J \to G/T$,
$\phi(x):=xT$ is surjective.
\end{la}
\begin{proof}
Let $d$ be the least common multiple
of~$t_\alpha$ (from Lemma~\ref{power})
and $d_J$.
Exploiting that $x^d \in D\cap M_J$
for every $x\in G$,
we can repeat the argument used to prove
Lemma~\ref{surjective}.
\end{proof}
The kernel of the homomorphism in  Lemma~\ref{MJnottrivial}
is equal to $D\cap M_J \cap T$ and so there is an exact sequence
%
%\ma{central}
\begin{equation}\label{central}
{\bf 1} \,\to \, D\cap M_J \cap T \, \to \,
D\cap M_J \, \to  \, G/T \, \to \, {\bf 1}\, .
\end{equation}
Recalling that $D$ is defined to be
$D = \wb{\langle x^{t_\alpha} \mid x\in G \rangle}$
and similarly for $\tilde D$,
the following lemma implies that
$q(D)\sub \tilde D$.
Hence $q(D\cap T)\sub \tilde D\cap \tilde T=\one$
and thus $D\cap T\sub S$,
whence $D\cap M_J\cap T$ is contained in the centre
of~$D\cap M_J$, by
Lemma~\ref{commute}. Thus, (\ref{central}) is
a central extension.
\begin{la}
$\tilde{t}_\alpha$ divides $t_\alpha$.
\end{la}
\begin{proof}
By the proof of Lemma~\ref{power},
$t_\alpha$ is the product
of the modules of the automorphisms
induced by $\alpha^{-1}$
on the composition factors of~$T$.
By Lemma~\ref{re-org},
the latter coincide with those composition factors
of~$G$ which are torsion groups.
Likewise,
$\tilde{t}_\alpha$
is the product
of the modules
of the automorphisms
induced by $\alpha^{-1}$
on those composition factors
of $G/S$ which are torsion groups.
As the latter are among
the composition factors
of $G$ which are torsion groups,
we deduce that~$\tilde{t}_\alpha$
divides~$t_\alpha$.
\end{proof}
Two algebraic results will help us to discuss
the central extension~(\ref{central}).
\begin{la}\label{tor-free-nilp}
Let $H$ be a
nilpotent group which admits a central series\linebreak
$\one=H_0\triangleleft H_1\triangleleft\cdots\triangleleft H_n=H$
such that $H/H_j$ is torsion-free
for each $j\in \{0,\ldots, n\}$.
Then roots in~$H$ are unique:
If $x^m=y^m$
for certain $x,y\in H$ and
$m\in \N$,
then $x=y$.
\end{la}
\begin{proof}
The proof is by induction on~$n$.
If $n=1$, then $H$ is abelian and torsion-free.
Thus $x^m=y^m$ entails
that $(y^{-1}x)^m=1$,
whence $y^{-1}x=1$ and $x=y$.
If $n>1$,
then the inductive hypothesis
applies to $H/H_1$,
whence $xH_1=yH_1$.
Since $H_1\sub Z(H)$,
we have $x=yz$
for some element $z$
in the centre of~$H$.
Thus $y^m=x^m=y^mz^m$
and hence $z^m=1$.
Since $H$ is torsion-free,
we infer that $z=1$ and thus $x=y$.
\end{proof}
%
%
%\ma{splits}
\begin{la}\label{splits}
Let ${\bf 1}\to C \to H \to H/C \to {\bf 1}$
be a central extension,
where $C$ is a group of finite exponent
and $Q:=H/C$ an infinitely divisible
nilpotent group admitting a central series
$\one=Q_0\triangleleft Q_1\triangleleft\cdots\triangleleft Q_m=Q$
such that $Q_j$ is infinitely divisible
and $Q/Q_j$ is torsion-free
for each $j\in \{0,\ldots, m\}$.
Then there is an infinitely divisible subgroup, $L$, of $H$
such that
\mbox{$H = C\times L$.}\linebreak
Furthermore, $C=\tor(H)$ and $L=\dv(H)$
are fully invariant subgroups~of~$H$.
\end{la}
\begin{proof}
It suffices to prove the first assertion
(because the final assertion is an
immediate consequence).
We first note that $C=\tor(H)$ because
$H/C$ is torsion-free.
Hence $C$ is a fully invariant subgroup.
Let $d\in \N$ be such that $c^d = 1$ for every $c\in C$.
We first show that each coset, $xC$, in $H/C$
contains a unique element that is divisible by $kd$
for every $k\in{\mathbb N}$.
Since $H/C$ is infinitely divisible, there is $yC\in H/C$
such that $(yC)^{kd} = xC$. Hence $y^{kd}\in xC$.
Because $H/C$ satisfies the hypotheses
of the preceding lemma,
the coset $yC$ is unique.
If $yc$ is another element of $yC$, then
$(yc)^{kd} = y^{kd}$ because $c$ belongs to the centre of $H$ and
$c^d = 1$. Hence $y^{kd}$ is the unique element of $xC$ that is divisible by
$kd$. If $k'$ is another element of ${\mathbb N}$, then there is $z\in
H$ such that $z^{kk'd} = y^{kd}$. Hence $y^{kd}$ is also the unique
element of $xC$ divisible by $k'd$.

Define $L := \{y^{d} \colon y\in H\}$.
By the above argument, each element of $L$
is divisible by $kd$ for each $k\in {\mathbb N}$.
Furthermore, $LC/C=H/C$.
We assert:
\emph{$L$~is a group and is complementary to $C$.}
The proof is by induction on~$m$.

Assume $m=1$ first; then $H/C$ is abelian.
Given $x_1,\,x_2\in L$,
we have
$x_i = y_i^{2d}$ for some $y_i\in H$.
Since $z := y_1^{-1}y_2^{-1}y_1y_2$ belongs to
$C$, we have $z^d= 1$
and thus
$x_1x_2 =
(y_1y_2)^{2d}
z^{\frac{1}{2}2d(2d-1)}
=
(y_1y_2)^{2d}
z^{d(2d-1)}
=(y_1y_2)^{2d}\in L$.
It is clear that $L$ is closed under inverses and so $L$ is a
group. It is also clear from its definition that $L$ is a
characteristic subgroup and that $L\cap C = {\bf 1}$.
Since, furthermore, $LC/C=H/C$,
we have $LC=H$ and hence $H = C\times L$.

Now let $m>1$ and assume that the assertion holds
if $m$ is replaced by $m-1$.
Let $H'$ be the inverse image of
$Q_1$ under the quotient map $H\to H/C$.
Then $H'/C=Q_1$
is an infinitely divisible, torsion-free
abelian group, whence the
extension ${\bf 1}\to C\to H'\to H'/C\to{\bf 1}$
satisfies the
hypotheses of the lemma.
By the abelian case already discussed,
$H' = C\times L'$
with $L'$ an infinitely divisible and characteristic subgroup of~$H'$.
The group $L'$ is normal in $H$ and the
extension $\one\to CL'/{L'}\to H/{L'}\to (H/L')/(CL'/L')\to\one$
satisfies the hypotheses of the lemma,
because $CL'/{L'}\isom C/(C\cap L')$
is a torsion group of finite exponent
and  $(H/L')/(CL'/L')\isom H/(CL')=H/{H'}\isom Q_m/Q_1$
is a nilpotent group
isomorphic to $Q_m/Q_1$ which
admits the central series
${\bf 1}=Q_1/Q_1\triangleleft\cdots \triangleleft Q_m/Q_1$
of length $m-1$ where $Q_j/Q_1$
is infinitely divisible
and $(Q_m/Q_1)/(Q_j/Q_1)\isom Q_m/Q_j$
torsion-free for all $j\in \{1,\ldots, m\}$.
Hence $H/L' =
(CL'/L')\times L''$ for an infinitely divisible subgroup~$L''$
of $H/{L'}$,
by the case $m-1$.
We claim:
\emph{The inverse image of $L''$ under the quotient map
$q\colon H\to H/L'$ is equal to~$L$.}
If this is true,
then $L$ is a characteristic subgroup
of~$H$. Because the elements of~$L$ are divisible
by $kd$ for each~$k$
and~$C=\tor(H)$ has finite exponent,
$L\cap C={\bf 1}$ holds.
Since, furthermore,
$LC/C=H/C$, we have $LC=H$ and hence
$H=C\times L$. In particular, $L\isom H/C$ is infinitely
divisible.

It only remains to verify the claim.
Since $CL'/L'$ has exponent
dividing~$d$, we have $q(y^d)\in L''$
for each $y\in H$ and thus $q(L)\sub L''$.
If $x\in L''$, then there exists
$w\in L''$ such that $w^d=x$.
Taking $y\in H$ such that
$w=q(y)$, we have $y^d\in L$ and
$q(y^d)=x$, showing that
$q(L)=L''$.
Hence $L=q^{-1}(L'')$
will follow if we can show that $LL'\sub L$.
To this end, let $x_1\in L$, $x_2\in L'$.
Then $x_1=y_1^{2d}$
for some $y_1\in H$, and $x_2=y_2^{2d}$
for some $y_2\in L'$.
Since $L'C/C\sub Q_1$ is contained
in the centre of $H/C$, we have
$z:=y_1^{-1}y_2^{-1}y_1y_2\in C$
and thus $L\ni (y_1y_2)^{2d}=y_1^{2d}y_2^{2d}z^{d(2d-1)}
=y_1^{2d}y_2^{2d}=x_1x_2$.
This completes the proof.
\end{proof} 
To obtain information concerning
the central extension~(\ref{central}),
we now consider the canonical $\alpha$-stable series
$\one=Q_0\triangleleft\cdots \triangleleft Q_m=G/T$
of $G/T$. Since all $\langle\alpha\rangle$-composition
factors of $G/T$ are torsion-free by Lemma~\ref{lemmins},
also all $\langle\alpha\rangle$-composition
factors of $Q_m/Q_j$ are torsion-free,
as well as those of~$Q_j$,
for $j\in \{0,\ldots, m\}$.
Hence $Q_j$ and $Q_m/Q_j$
are torsion-free, infinitely divisible,
nilpotent groups by Proposition~\ref{struc-t-free}
and Lemma~\ref{div-torfree}.
Thus
Lemma~\ref{splits} implies that
(\ref{central}) splits as an extension
of abstract groups. Write
\[
D\cap M_J \; = \; (D\cap M_J \cap T)\times L,
\] 
say, where $L \cong G/T$.
Then $L=\dv(D\cap M_J)$
is a characteristic subgroup
of $D\cap M_J$. Since $D\cap M_J$
is a normal and $\alpha$-stable
subgroup of $G$, it
follows that $L$ is a normal and $\alpha$-stable
subgroup of $G$. We also have that $G
= LT$, by Lemma~\ref{MJnottrivial}, and that $L\cap T = {\bf 1}$,
since~$L$ is torsion-free.
Therefore, $G = T\times L$
as an abstract group. We need more:
%
%
%\ma{newlm}
\begin{la}\label{newlm}
$L$ is closed in~$G$
and $G=T\times L$ as a topological group.
\end{la}
\begin{proof}
Pick a compact, open
subgroup $U\leq G$
and let $\pi\colon G\to G/T$ be the quotient map.
Since $G/T$ is a product
of $p$-adic Lie groups for certain primes~$p$,
we find primes
$p_1,\ldots,p_m$ and continuous homomorphisms
$\xi_i\colon \Z_{p_i}\to G/T$
such that, for each $k\in \N_0$,
\[
V_k\; :=\;
\xi_1(p_1^k\Z_{p_1})\,
\xi_2(p_2^k\Z_{p_2})\, \cdots
\, \xi_m(p_m^k\Z_{p_m})
\]
is a compact, open subgroup of~$\pi(U)$.
To see this, recall that
each $p$-adic Lie group
admits coordinates of the second kind,
and apply the ultrametric inverse function theorem.
Since $\pi|_U\colon U\to \pi(U)$
is a quotient homomorphism between
pro-finite groups, each $\xi_i$ lifts
to a continuous homomorphism
$\theta_i\colon \Z_{p_i}\to U$
such that $\pi\circ\theta_i=\xi_i$
(this is clear from
standard facts of
pro-finite Sylow theory,
notably \cite[Proposition~2.2.3\,(b)]{Wsn}).
There exists $k\in \N$
such that $t_\alpha^{-1}p_i^k\in \Z_{p_i}$
for all $i\in \{1,\ldots, m\}$,
entailing that all elements of
$\theta_i(p^k_i\Z_{p_i})$
are divisible by~$t_\alpha$ in~$G$.
Hence
$\theta_i(p^k_i\Z_{p_i})\sub L$
and thus
\[
W\; :=\;
\theta_1(p_1^k\Z_{p_1})\,
\theta_2(p_2^k\Z_{p_2})\, \cdots\,
\theta_m(p_m^k\Z_{p_m})\; \sub\; L\,.
\]
Note that $\pi(W)=V_k$.
If $x,y\in W$,
then $\pi(x)\pi(y)\in V_k$
and thus $\pi(x)\pi(y)=\pi(z)$
for some $z\in W$.
Since also $\pi(xy)=\pi(z)$,
where both $xy$ and $z$ are in~$L$,
we deduce from the injectivity of
$\pi|_L$ that $xy=z$.
Hence $W$ is a subgroup of~$L$.
Since $W$ is compact, the bijective continuous
homomorphism $\pi|_W^{V_k}$
is an isomorphism of topological
groups.
Therefore the homomorphism
$(\pi|_L)^{-1}\colon G/T\to L$
is continuous on the identity neighbourhood~$V_k$
and hence continuous.
Thus $G=T\times L$ as a topological
group, and thus $L$~is~closed.
\end{proof}
{\bf Proof of Lemma~\ref{factor}, completed.}
Since $G=T\times L$,
Lemma~\ref{power}
implies that $x^{t_\alpha} \in L$ for each
$x\in G$.
Thus $D\stackrel{\text{def}}{=}\wb{\langle x^{t_\alpha}\colon x\in G\rangle} \sub \wb{L}=L$,
using Lemma~\ref{newlm}.
Conversely, $L\sub D$,
because~$L$ is infinitely divisible.
Hence $L=D$ and thus $G =
T\times D$ internally as a topological group.
Since $T$ has finite exponent
and $D$ is infinitely divisible,
it follows that $D=\dv(G)$.\Punkt
\appendix
\section{Appendix: Proof of Lemma~\ref{abstrgp}.}\label{appdetails}
%
%\ma{appdetails}
%
The case $n=1$ being trivial
(defining the empty intersection as~$G$),
we may assume that $n>1$ and that
the claim holds for $n-1$
distinct normal subgroups.
It is clear that $\theta$ is
injective and takes $D_k/D$ into $G/N_k$,
for each~$k$; it therefore
only remains to show that $\theta(D_k/D)=G/N_k$,
i.e., $D_kN_k/N_k=G/N_k$.
If this is wrong, then
$D_kN_k/N_k$ is a proper normal
subgroup of the simple group $G/N_k$
for some $k\in \{1,\ldots, n\}$,
whence $D_kN_k=N_k$
and thus $D_k\sub N_k$.
Then $N_k/D_k$ is a proper normal subgroup
of $G/D_k$. By induction,
$G/D_k$ is the internal direct
product of the non-abelian simple
groups
$C_j/D_k\isom (G/D_k)/(N_j/D_k) \isom G/N_j$
for $j\in \{1,\ldots, n\}\setminus\{k\}$,
where $C_j$ is the intersection of
the groups $N_i$ for all
$i\in \{1,\ldots, n\}\setminus\{j,k\}$.
By Remak's Theorem \cite[3.3.12]{Rob},
there exists a finite subset $F\sub \{1,\ldots, n\}\setminus
\{k\}$ such that $N_k/D_k=\prod_{j\in F}(C_j/D_k)$
internally. Since $N_k/D_k$ is a proper
subgroup of $G/D_k$,
there exists $\ell\in \{1,\ldots, n\}\setminus\{k\}$
such that $\ell\not\in F$.
Then $C_j\sub N_\ell$ for each $j\in F$
and thus $N_k/D_k\sub N_\ell/D_k$,
entailing that $N_k\sub N_\ell$.
But then $N_\ell/N_k$ is a proper normal
subgroup of the simple group $G/N_k$
and thus $N_\ell/N_k=\one$.
Hence $N_\ell=N_k$ and therefore
$\ell=k$, contradicting our choice
of~$\ell$. This completes the
inductive proof.\vspace{2.5mm}\Punkt

{\footnotesize
{\bf Helge Gl\"{o}ckner}, TU Darmstadt, Fachbereich Mathematik AG~5,
Schlossgartenstr.\,7,\\
64289 Darmstadt, Germany. E-Mail:
gloeckner\at{}mathematik.tu-darmstadt.de\\[2mm]
{\bf George A. Willis},
Department of Mathematics,
University of Newcastle,
Callaghan,\\
NSW 2308, Australia.
E-Mail:
George.Willis\at{}newcastle.edu.au}
\end{document}